\documentclass{amsart}

\newtheorem{thm}{Theorem}[section]
\newtheorem{prop}[thm]{Proposition}
\newtheorem{cor}[thm]{Corollary}
\newtheorem{lem}[thm]{Lemma}

\theoremstyle{remark}
\newtheorem{rem}[thm]{Remark}

\numberwithin{equation}{section}

\def \vtk{V^{\otimes k}}
\def \vti{V^{\otimes k_i}}
\def \bY{\bar Y}
\def \D{\Delta_k}
\def \Z{\mathbb Z}
\def \C{\mathbb C}
\def \Q{\mathbb Q}
\def \N{\mathbb N}
\def \wt{{\rm wt}}
\def \Res{{\rm Res}}
\def \End{{\rm End}}
\def \Aut{{\rm Aut}}
\def \mod{{\rm mod}}
\def \o{\omega}
\def \l{\lambda}

\newcommand{\ZZ}[0]{{\Z}_+}
\newcommand{\Lx}[0]{x^{j + 1} \frac{\partial}{\partial x}}
\newcommand{\Lo}[0]{x \frac{\partial}{\partial
x}}  
\newcommand{\Ly}[0]{y^{j + 1} \frac{\partial}{\partial y}}
\newcommand{\Loy}[0]{y \frac{\partial}{\partial
y}} 

\begin{document}

\title[Permutation twisted tensor product VOAs]{Twisted sectors for tensor
product vertex operator algebras associated to permutation groups}

\author{Katrina Barron}
\address{Department of Mathematics, University of California, Santa Cruz, CA 95064}
\email{kbarron@math.ucsc.edu}
\thanks{The first author was supported by an NSF Mathematical Sciences Postdoctoral
Research Fellowship and a University of California President's Postdoctoral Fellowship.}

\author{Chongying Dong}
\address{Department of Mathematics, University of California, Santa Cruz, CA 95064}
\email{dong@math.ucsc.edu}
\thanks{The second author was supported by NSF grant DMS-9700923 and a research grant {}from 
the Committee on Research, UC Santa Cruz.}

\author{Geoffrey Mason}
\address{Department of Mathematics, University of California, Santa Cruz, CA 95064}
\email{gem@math.ucsc.edu}
\thanks{The third author was supported by NSF grant DMS-9700909 and a research
 grant {}from the Committee on Research, UC Santa Cruz.}




\begin{abstract}
Let $V$ be a  vertex operator algebra, and for $k$ a 
positive integer, let $g$ be a $k$-cycle permutation of the vertex
operator algebra $V^{\otimes k}$.  We prove that the categories of
weak, weak admissible and ordinary $g$-twisted modules for the tensor
product vertex operator algebra $V^{\otimes k}$ are isomorphic to the 
categories of weak, weak admissible and ordinary $V$-modules, 
respectively.  The main result is an explicit construction of the 
weak $g$-twisted $V^{\otimes k}$-modules {}from weak $V$-modules.  For 
an arbitrary permutation automorphism $g$ of $V^{\otimes k}$ the 
category of weak admissible $g$-twisted modules for $V^{\otimes k}$ 
is semisimple and the simple objects are determined if $V$ is 
rational.  In addition, we extend these results to the more general 
setting of $\gamma g$-twisted $V^{\otimes k}$-modules for $\gamma$ a 
general automorphism of $V$ acting diagonally on $V^{\otimes k}$ and 
a $g$ a permutation automorphism of $V^{\otimes k}$.
\end{abstract}

\maketitle

\section{Introduction}

Orbifold theory and coset construction theory \cite{GKO} are two 
important ways of constructing a new conformal field theory {}from a 
given one. The first orbifold conformal field theory was introduced 
in \cite{FLM1} and the theory of orbifold conformal fields was 
subsequently developed, for example, in \cite{DHVW1}, \cite{DHVW2}, 
\cite{FLM3}, and \cite{DVVV}.  The first critical step in orbifold 
conformal field theory is to construct twisted sectors.  In 
\cite{Le1} and \cite{FLM2}, twisted sectors for finite automorphisms 
of even lattice vertex operator algebras were first constructed -- 
the twisted vertex operators were constructed, and in \cite{Le2} and 
\cite{DL2}, the twisted Jacobi identity was formulated and shown to 
hold for these operators.  These results led to the introduction of
the notion of $g$-twisted $V$-module \cite{FFR}, \cite{D}, for $V$ a 
vertex operator algebra and $g$ an automorphism of $V$.  This notion 
records the properties obtained in \cite{Le1}, Section 3.3 of 
\cite{FLM2}, and \cite{Le2}, and provides an axiomatic definition of 
twisted sectors.

The focus of this paper is the construction of twisted sectors for 
permutation orbifold theory.  Let $V$ be a  vertex operator algebra, 
and for a fixed positive integer $k$,  consider the tensor product 
vertex operator algebra $\vtk$ (see \cite{FLM3}, \cite{FHL}).  Any 
element $g$ of the symmetric group $S_k$ acts on $\vtk$ in the 
obvious way, and thus it is appropriate  to consider $g$-twisted 
$\vtk$-modules.  This is the setting for permutation orbifolds.  In 
the case of $V$ a lattice vertex operator algebra, this becomes a 
special case of the more general result of \cite{Le1}, \cite{Le2}, 
\cite{FLM2}, and \cite{DL2}, mentioned above.   Permutation 
orbifold theory has also been studied in the physics  literature, 
e.g., in \cite{KS} and \cite{FKS}, and in \cite{BHS}, the twisted 
vertex operators for the generators were given for affine, Virasoro, 
super-Virasoro and $W_3$ cyclic permutation orbifolds.  The  
characters and modular properties of permutation orbifolds are  
presented in \cite{Ba}.  However, the construction of twisted sectors 
for general permutation orbifold theory and for arbitrary orbifold 
theory have been open problems.  The main result of this paper is the 
explicit construction of twisted sectors for general permutation 
orbifold theory.  In addition, for $g$ a $k$-cycle, we show that the 
categories of weak, weak admissible and ordinary $g$-twisted
$\vtk$-modules are isomorphic to the categories of weak, weak 
admissible and ordinary $V$-modules, respectively. (The definitions 
of weak, weak admissible and ordinary twisted modules are given in  
Section 3.)  Our proof includes an explicit construction of the 
$g$-twisted $\vtk$-modules given a $V$-module.  We show that for an 
arbitrary permutation automorphism $g$ of $V^{\otimes k}$ the 
category of weak admissible $g$-twisted modules for $V^{\otimes k}$ 
is semisimple and the simple objects are determined if $V$ is rational.  
In addition, we extend our results to the more general setting of 
$\gamma g$-twisted $V^{\otimes k}$-modules for $\gamma$ a general 
automorphism of $V$ acting diagonally on $V^{\otimes k}$ and $g$ an 
arbitrary permutation automorphism of $V^{\otimes k}$.  One can use  
our constructions to calculate the characters and perform modular  
transformations.  Moreover, we expect that the methods introduced in 
this paper can be extended to construct twisted sectors for arbitrary 
orbifold theory.

We next point out some recent results and conjectures in the theory
of vertex operator algebras with which we hope to put the results of 
this paper into perspective and show some of the motivation which led 
to our results.  Let $V$ be a vertex operator algebra and $g$ an 
automorphism of $V$.  In \cite{DM} it is shown that given a weak 
$V$-module $(M,Y)$, one can define a new weak $V$-module
$g\circ M$ such that $g\circ M=M$ as the underlying space and
the vertex operator associated to $v$ is given by $Y(gv,z)$.  Then $M$
is called {\it $g$-stable} if $g\circ M$ and $M$ are isomorphic as
weak $V$-modules.  It is a well-known {\em conjecture} that if $V$ is
rational and $g$ is of finite order, then the number of isomorphism
classes of irreducible $g$-stable $V$-modules is finite and equal to
the number of isomorphism classes of irreducible $g$-twisted
$V$-modules.  It is proved in \cite{DLM3} that if $V$ is rational and
satisfies the $C_2$ condition, and $g$ is of finite order, then the
number of isomorphism classes of irreducible $g$-twisted modules is
finite and less than or equal to the number of isomorphism classes of
irreducible $g$-stable $V$-modules. Moreover, if $V$ is also assumed
to be $g$-rational, then the number of isomorphism classes of
irreducible $g$-twisted $V$-modules is equal to the number of 
isomorphism classes of irreducible $g$-stable $V$-modules.

Now consider the tensor product vertex operator algebra $\vtk$ as
discussed above with $g$ a $k$-cycle. {}From Proposition 4.7.2 and
Theorem 4.7.4 in \cite{FHL}, it follows that the number of isomorphism
classes of irreducible $g$-stable $\vtk$-modules is equal to the
number of isomorphism classes of irreducible $V$-modules.  Combining
this result with the conjecture above, we have the conjecture that if
$V$ is rational and $g$ is a $k$-cycle, then the number of isomorphism
classes of irreducible $g$-twisted $\vtk$-modules is equal to the
number of isomorphism classes of irreducible $V$-modules. Thus we see
that the main result in this paper is actually a stronger result than
this conjecture, namely that the categories of weak, weak admissible 
and ordinary $g$-twisted $\vtk$-modules are in fact {\it isomorphic} 
to the categories of weak, weak admissible and ordinary $V$-modules,
respectively, even without the assumption that $V$ is rational.

To construct the isomorphism between the category of weak 
$g$-twisted $\vtk$-modules and the category of weak  
$V$-modules for $g$ a $k$-cycle, we define a weak  
$g$-twisted $\vtk$-module structure on any weak $V$-module.
Our construction has been motivated by the following two results.

The first motivating result is the modular invariance of the trace
function in orbifold theory \cite{DLM3}.  If $V$ is a holomorphic vertex
operator algebra (see \cite{DLM3}), then $V$ is the only irreducible module
for itself.  Thus by Proposition 4.7.2 and Theorem 4.7.4 in \cite{FHL}, 
$\vtk$ is the only irreducible $g$-stable $\vtk$-module. Following the 
notation and results of \cite{DLM3}, consider the symmetric group $S_k$ as 
an automorphism group of $\vtk$ and denote by $Z(x,y,\tau)$ the $y$-trace 
on the unique $x$-twisted $\vtk$-module $\vtk(x)$ where $x,y\in S_k$
commute. Then the span of $Z(x,y,\tau)$ is modular invariant.  Using
the modular invariance result of \cite{DLM3}, one can show that the
graded dimension of irreducible $g$-twisted $\vtk$-modules is exactly
the graded dimension of $V$ except that $V_n$ is graded by $n/k$ (plus
a uniform shift) instead of $n$ for $n\in\Z$. (This fact has also been
observed and used in \cite{DMVV} to study elliptic genera of symmetric
products and second quantized strings.)  This leads one to expect a
$g$-twisted $\vtk$-module structure on $V$ with the new gradation.

The second motivating result is the construction in \cite{Li2} of certain
$g$-twisted $V$-modules for $g$ a certain automorphism of $V$.  Let
$h\in V_1$ such that the zero-mode operator $h_0$ for $Y(h,z)$ acts
semisimply on $V$ with only finitely many eigenvalues and such that
these eigenvalues are rational.  Then $e^{2\pi h_0 i}$ is an
automorphism of $V$.  Li's construction defines a new action on any
$V$-module and this action gives a $e^{2\pi h_0 i}$-twisted $V$-module
structure on the original $V$-module. The main feature in the new
action is an exponential operator $\Delta(z)$ built up {}from the
component operators $h_n$ of $Y(h,z)$ for $n\geq 0$.  This kind of
operator first appeared in the construction of twisted sectors for 
lattice vertex operator algebras \cite{FLM2}, \cite{FLM3}.

Now let $M$ be a weak admissible $V$-module and let $g$ be the
$k$-cycle $g = (12 \cdots k)$.  As suggested above {}from our first
motivating result, one expects a weak $g$-twisted $\vtk$-module 
structure on $M$ with some vertex operator $Y_g(v,z)$ for $v\in\vtk$. 
{}From the twisted Jacobi identity one sees that the component 
operators of $Y_g(u\otimes {\bf 1}\otimes \cdots\otimes {\bf 1},z)$ 
for $u\in V$ form a Lie algebra.  In addition, the fact that 
$Y_g(gv,z) = \lim_{z^{1/k}\to \eta^{-1}z^{1/k}}Y_g(v,z)$, for $\eta 
= e^{-2 \pi i/k}$, indicates that all the  vertex operators 
$Y_g(v,z)$ for $v\in \vtk$ are generated by $Y_g(u\otimes {\bf 1}
\otimes \cdots\otimes {\bf 1},z)$ for $u\in V$. Therefore, the key 
point is to define $Y_g(u\otimes {\bf 1}\otimes \cdots\otimes {\bf 1},
z)$ which is expected to be $Y(\D(z)u,z^{1/k})$ for some operator 
$\D(z)\in (\End \; V)[[z^{1/k},z^{-1/k}]]$ due to our second 
motivating result above.  In \cite{FLM2}, \cite{FLM3} and \cite{Li2} 
an operator $\Delta(z)$ was introduced by using vertex operators 
associated to certain Heisenberg algebras.  But for an arbitrary 
vertex operator algebra, we do not have a Heisenberg algebra 
available.  So we must find another way to construct $\Delta(z)$. 
Note that such a construction should also work if $V$ is the vertex 
operator algebra associated to the highest weight modules for the 
Virasoro algebra.  Therefore in general, the only available operators 
with which we can construct $\Delta(z)$ are $L(n)$ for $n\in \Z$.  And 
in fact, our operator $\D(z)$ is built up {}from $L(n)$ for $n\geq 0$ 
(see Section 2).

The paper is organized as follows. In Section 2, we define the
operator $\D(z)$ on $V$ and prove several important properties of
$\D(z)$ which are needed in subsequent sections.  The main ideas for the
proofs of these identities come {}from \cite{H}. In Section 3, we define a
weak $g$-twisted $\vtk$-module structure on any weak $V$-module $M$
for $g = (12 \cdots k)$ by using the operator $\D(z)$.  As a result we
construct a functor $T_g^k$ {}from the category of weak $V$-modules to
the category of weak $g$-twisted $\vtk$-modules such that $T_g^k$ maps
weak admissible (resp., ordinary) $V$-modules into weak admissible
(resp., ordinary) $g$-twisted $\vtk$-modules.  In addition, $T_g^k$
preserves irreducible objects. In Section 4, we define a weak
$V$-module structure on any weak $g$-twisted $\vtk$-module.  In so
doing, we construct a functor $U_g^k$ {}from the category of weak
$g$-twisted $\vtk$-modules to the category of weak $V$-modules such
that $T_g^k \circ U_g^k = id$ and $U_g^k \circ T_g^k = id$.  In Section
5, we give the extension of the results of Section 3 and 4 for $g$ a
general $k$-cycle.

Section 6 is devoted to twisted modules for an arbitrary permutation 
$g \in S_k$.  In particular we prove that if $V$ is rational, then 
$\vtk$ is $g$-rational.  We also construct irreducible $g$-twisted 
$\vtk$-modules {}from irreducible $V$-modules.

In Section 7, we study various twisted modules for an automorphism
of $\vtk$ which is a product of a permutation and an automorphism of
$V$.  Here the automorphisms of $V$ act on $\vtk$ diagonally. For 
$\gamma$ an automorphism of $V$, $g$ the $k$-cycle $g =(1 2 \cdots k)$, 
and $\gamma g$ an automorphism of $\vtk$, we  show that the category 
of weak $\gamma g$-twisted $\vtk$-modules is isomorphic to the category 
of weak $\gamma^k$-twisted $V$-modules and that this isomorphism 
preserves admissible, ordinary and irreducible objects. Finally we
construct $\gamma g$-twisted $\vtk$-modules {}from $\gamma$-twisted
$V$-modules for $g$ an arbitrary permutation and show that if $V$ is 
$\gamma$-rational, then $\vtk$ is $\gamma g$-rational.

The authors thank James Lepowsky for pointing out some mistakes in 
earlier versions of this paper and for giving helpful comments on the
paper's exposition.  We would also like to thank Hirotaka Tamanoi for
valuable discussions.

\section{The operator $\Delta_k (z)$}
\setcounter{equation}{0} 

In this section we define an operator $\Delta_k(z)=\Delta_k^V(z)$ on
a vertex operator algebra $V$ for a fixed positive integer $k$. In 
Section 3, we will use $\Delta_k(z)$ to construct a $g$-twisted 
$V^{\otimes k}$-module {}from a $V$-module where $g$ is a certain $k$-cycle.

Let $\ZZ$ denote the positive integers.  Let $x$, $y$, $z$, $z_0$, and 
$\alpha_j$ for $j \in \ZZ$ be formal variables commuting with each other.  
Consider the polynomial
\[\frac{1}{k} (1 + x)^k - \frac{1}{k} \in x \C [x] . \]  
By Proposition 2.1.1 in \cite{H}, for any formal power series $\sum_{j
\in \ZZ} c_j x^j \in x \C [[x]]$ there exist unique $a_j \in \C$
for $j \in \ZZ$ such that
\[\exp \Biggl( - \sum_{j \in \ZZ} a_j  \Lx
\Biggr) \cdot x = \sum_{j \in \ZZ} c_j x^j .\]
Thus for $k \in \ZZ$, we can define $a_j \in \C$ for $j \in \ZZ$, by  
\[\exp \Biggl( - \sum_{j \in \ZZ} a_j  \Lx
\Biggr) \cdot x = \frac{1}{k} (1 + x)^k -
\frac{1}{k} .\]
For example, $a_1=(1-k)/2$ and $a_2=(k^2-1)/12.$

Let $R$ be a ring, and let $O$ be an invertible linear operator on
$R[x, x^{-1}]$.  We define another linear operator $O^{\Lo}$ by
\[O^{\Lo} \cdot x^n = O^n x^n \]
for any $n \in \Z$.  For example, since $z^{1/k}$ 
can be thought of as an invertible linear multiplication operator  
$\C [x, x^{-1}]$, we have the operator $z^{(1/k) \Lo}$ {}from $\C 
[x,x^{-1}]$ to $\mathbb{C}[z^{1/k},z^{-1/k}]
[x,x^{-1}]$.  Note that $z^{(1/k) \Lo}$ can  be extended to a 
linear operator on $\C [[x,x^{-1}]]$ in the obvious way.
  
Let
\begin{eqnarray*}
f(x) &=& z^{1/k} \exp \Biggl(- \sum_{j \in \ZZ} a_j  \Lx \Biggr) 
\cdot x \\
&=& \exp \Biggl(- \sum_{j \in \ZZ} a_j  \Lx \Biggr) \cdot
z^{(1/k) \Lo} \cdot x\\
&=& \frac{z^{1/k}}{k} (1 + x)^k - \frac{z^{1/k}}{k} \; \; \in
z^{1/k}x\mathbb{C}[x].
\end{eqnarray*}
Then the compositional inverse of $f(x)$ in $x\mathbb{C}[z^{-1/k}, 
z^{1/k}][[x]]$ is given by
\begin{eqnarray*}
f^{-1} (x) &=& z^{- (1/k) \Lo} \exp \Biggl( \sum_{j \in \ZZ} a_j  \Lx 
\Biggr) \cdot x \\
&=& z^{- 1/k} \exp \Biggl( \sum_{j \in \ZZ} a_j 
z^{- j/k} \Lx \Biggr) \cdot x \\
&=& (1 + k z^{- 1/k} x)^{1/k} - 1
\end{eqnarray*} 
where the last line is considered as a formal power series in
$z^{-1/k}x \mathbb{C}[z^{-1/k}][[x]]$, i.e., we are expanding
about $x = 0$ taking $1^{1/k} = 1$.

Let $f_\alpha (x)$ denote the formal power series 
\[f_\alpha (x) = z^{1/k} \exp \Biggl(- \sum_{j \in \ZZ} \alpha_j  
\Lx \Biggr) \cdot x \in z^{1/k}x + z^{1/k}x^2 
\C [\alpha_1, \alpha_2,...]  [[x]] .\] 
Then 
\[f^{-1}_\alpha (x) \in z^{-1/k}x + x^2 \C [x]
[z^{-1/k}][[\alpha_1, \alpha_2,...]],\]
and 
\[f_\alpha
(f^{-1}_\alpha (x) + z^{-1/k} y) - x \in \C[x]
[z^{1/k}, z^{-1/k}][[\alpha_1, \alpha_2, ...]][[y]].\]
Furthermore, the coefficient of the monomial $y$ in $f_\alpha
(f^{-1}_\alpha (x) + z^{-1/k} y) - x$ is in $1 + x^2 \C[x]
[z^{1/k}, z^{-1/k}][[\alpha_1, \alpha_2, ...]]$.  
Therefore, following \cite{H}, we can define 
\[\Theta _j = \Theta_j (\{- \alpha_n \}_{n \in \ZZ},
z^{1/k}, x) \in \C[x][z^{1/k}, 
z^{-1/k}][[\alpha_1, \alpha_2,...]]\] 
for $j \in {\N}$ by
\[ e^{\Theta_0}\exp \Biggl( \sum_{j \in \ZZ} \Theta_j \Ly \Biggr)  y 
= f_\alpha (f^{-1}_\alpha (x) + z^{-1/k} y) - x .\]

\begin{prop}\label{Theta prop}
$\Theta_j (\{- a_n \}_{n \in \ZZ}, z^{1/k}, \frac{1}{k}
z^{1/k - 1} z_0)$ for $j \in {\N}$ is well defined in
$\C[z_0][[z^{-1/k}]]$.  Furthermore
\begin{equation}\label{Theta j}
\Theta_j (\{- a_n \}_{n \in \ZZ}, z^{1/k}, \frac{1}{k}
z^{1/k - 1} z_0) = - a_j (z + z_0)^{- j/k} 
\end{equation} 
for $j \in \ZZ$, and
\begin{equation}\label{Theta zero}
\exp \left(\Theta_0 (\{- a_n \}_{n \in \ZZ}, z^{1/k}, \frac{1}{k}
z^{1/k - 1} z_0) \right) = z^{ 1/k - 1} (z +
z_0)^{- 1/k + 1} , 
\end{equation} 
where $(z + z_0)^{- j/k}$ is understood to be expanded in 
nonnegative integral powers of $z_0$.
\end{prop}

\begin{proof} By Lemma 4.3.4 in \cite{H}, the formal series
$\Theta_j (\{- \alpha_n \}_{n \in \ZZ}, z^{1/k}, x)$ for $j
\in {\N}$, are actually in $\C[x][\alpha_1, \alpha_2, ...]  [[z^{-
1/k} ]]$.  Therefore $\Theta_j (\{- a_n \}_{n \in \ZZ},
z^{1/k}, x)$ is well defined in $\C[x][[z^{- 1/k}
]]$, for $j \in \mathbb{N}$, and the first statement of the proposition 
follows.

In $y \C[z_0] [[z^{-1/k}]][[y]]$, we have 
\begin{eqnarray*} 
& &\hspace{-.4in} z^{1/k - 1} (z + z_0)^{-1/k + 1}
\exp \Biggl(- \! \sum_{j \in \ZZ} a_j (z + z_0)^{-j/k} \Ly
\Biggr) \cdot y =  \\
&=& \! z^{1/k - 1} (z + z_0)^{-1/k + 1} 
(z + z_0)^{-(1/k) \Loy} (z + z_0)^{1/k}
\exp \Biggl(- \! \sum_{j \in \ZZ} a_j \Ly \Biggr) \! \cdot y \\
&=& \! z^{-1} (z + z_0) (z + z_0)^{-(1/k) \Loy} f(y)\\
&=& \! z^{-1} (z + z_0)  f((z + z_0)^{-1/k} y)\\
&=& \! z^{-1} (z + z_0) \Bigl( \frac{z^{1/k}}{k}  
(1 + (z + z_0)^{-1/k} y)^k - \frac{z^{1/k}}{k}  
\Bigr) \\
&=& \! \left. \frac{z^{1/k}}{k} \left(z^{-1/k} y + 
z^{-1/k} (z + k z^{1 - 1/k} x)^{1/k} \right)^k 
- \frac{ z^{1/k}}{k}  - x \right|_{x =\frac{1}{k} 
z^{1/k - 1} z_0}\\
&=& \! \left. f(z^{-1/k} (z + k z^{1 - 1/k} 
x)^{1/k} - 1 + z^{-1/k}y ) - x \right|_{x =\frac{1}{k} 
z^{1/k - 1} z_0}\\
&=& \! \left. f(f^{-1} (x) + z^{-1/k}y ) - x \right|_{x =\frac{1}{k} 
z^{1/k - 1} z_0}\\ 
&=& \! \exp \left(\Theta_0 (\{- a_n \}_{n \in \ZZ}, z^{1/k}, 
x) \right) \cdot \\
& & \! \hspace{1in} \left. \exp \Biggl( \sum_{j \in \ZZ} \Theta_j 
(\{- a_n \}_{n \in \ZZ}, z^{1/k}, x) \Ly \Biggr) \cdot y 
\right|_{x =\frac{1}{k} z^{1/k - 1} z_0} .
\end{eqnarray*}
Equations (\ref{Theta j}) and (\ref{Theta zero}) follow. \end{proof}

Let $V= (V, Y, {\bf 1}, \omega)$ be a vertex operator algebra.  In 
$(\End \;V)[[z^{1/k}, z^{-1/k}]]$, define 
\[\Delta_k^V (z) = \exp \Biggl( \sum_{j \in \ZZ} a_j z^{- j/k} L(j) 
\Biggr) k^{-L(0)} z^{\left( 1/k - 1 \right) L(0)} .\] 

\begin{prop}\label{psun1} 
In $(\End \;V)[[z^{1/k}, z^{-1/k}]]$, we have
\[\Delta_k^V (z) Y(u, z_0) \Delta_k^V (z)^{-1} = Y(\Delta_k^V 
(z + z_0)u, \left( z + z_0 \right)^{1/k} 
- z^{1/k}  ) ,\]
for all $u \in V$.
\end{prop}

\begin{proof} By equation (5.4.10)\footnote{There is a typo in this equation in
\cite{H}.  $A^{(0)}$ in the first line of equation (5.4.10) should be
$A^{(1)}$ which is the infinite series $\{A_j^{(1)}\}_{j\in\Z_{+}}$,
where $A_j^{(1)}\in \C$. In our case, $A_j^{(1)}=-a_j.$} in \cite{H}
and Proposition \ref{Theta prop} above, we have
\begin{eqnarray*}
& & \hspace{-.4in} \Delta_k^V (z) Y(u, z_0) \Delta_k^V (z)^{-1} = \\
&=& \! \exp \Biggl(\sum_{j \in \ZZ} a_j z^{- j/k} L(j) \Biggr) 
Y(k^{- L(0)} z^{\left( 1/k - 1 \right) L(0)} u, 
\frac{1}{k} z^{ 1/k - 1}z_0)\cdot \\
& & \hspace{3in} \cdot\exp \Biggl(- \sum_{j \in \ZZ} 
a_j z^{- j/k} L(j) \Biggr) \\ 
&=& \! z^{1/kL(0)} \exp \Biggl(\sum_{j \in \ZZ} 
a_j L(j) \Biggr) z^{-(1/k)L(0)} 
Y(k^{- L(0)} z^{\left( 1/k - 1 \right) L(0)} u, \frac{1}{k} 
z^{ 1/k - 1}z_0)\cdot \\
& & \hspace{2.1in} \cdot z^{(1/k)L(0)} \exp \Biggl( - 
\sum_{j \in \ZZ} a_j L(j) \Biggr) z^{-(1/k)L(0)} \\
&=& \! z^{(1/k)L(0)} Y \Biggl( z^{-(1/k)L(0)} 
\exp \Biggl( - \sum_{j \in \ZZ} \Theta_j(\{- a_n \}_{n \in \ZZ}, 
z^{1/k},\frac{1}{k} z^{ 1/k - 1}z_0  ) L(j) \Biggr) 
\cdot \Biggr. \\
& & \hspace{.2in} \cdot\exp \left(- \Theta_0 (\{- a_n \}_{n \in \ZZ}, 
z^{1/k}, \frac{1}{k} z^{ 1/k - 1}z_0) L(0) \right) 
\cdot k^{- L(0)} z^{\left( 1/k - 1 \right) L(0)} u, \\
& & \hspace{2.9in} \Biggl. f^{-1} (\frac{1}{k} z^{ 1/k - 1} z_0 ) 
\Biggr) z^{-(1/k)L(0)}\\
&=& \! z^{(1/k)L(0)} Y \Biggl(z^{-(1/k)L(0)} \exp 
\Biggl( \sum_{j \in \ZZ} a_j (z + z_0)^{- j/k} L(j) \Biggr) 
z^{-\left(1/k - 1 \right) L(0)} \cdot \Biggr. \\
& & \hspace{0.2in} \Biggl. \cdot 
(z + z_0)^{\left(1/k - 1 \right) L(0)} \cdot k^{- L(0)} 
z^{\left( 1/k - 1 \right) L(0)} u, \left( 1 + z^{- 1}z_0 \right)^{1/k} - 1 
\Biggr) z^{-(1/k)L(0)} \\
&=& \! z^{(1/k)L(0)} Y \left(z^{-(1/k)L(0)} \Delta_k^V (z + z_0) u, 
\left( 1 + z^{- 1}z_0 \right)^{1/k} - 1 \right) 
z^{-(1/k)L(0)}\\
&=& \! Y \left(\Delta_k^V (z + z_0) u, 
\left( z + z_0 \right)^{1/k} - z^{1/k}  \right) 
\end{eqnarray*}
as desired. \end{proof}

Define $\Delta_k^x (z) \in (\End \; \C[x, x^{-1}])
[[z^{1/k}, z^{-1/k}]]$ by
\[\Delta_k^x (z) = \exp \Biggl( - \sum_{j \in \ZZ} a_j z^{- j/k} \Lx
\Biggr) k^{\Lo} z^{\left( - 1/k + 1 \right) \Lo} . \]

\begin{prop}
In $(\End \; \C[x, x^{-1}]) [[z^{1/k}, 
z^{-1/k}]]$, we have 
\begin{eqnarray}
- \Delta_k^x (z) \frac{\partial}{\partial x} + \frac{1}{k} 
z^{1/k - 1} \frac{\partial}{\partial x}
\Delta_k^x (z)  &=& \frac{\partial}{\partial z} \Delta_k^x (z) , 
\label{identity in rep}\\
- \Delta_k^x (z)^{-1} \frac{\partial}{\partial x} + k 
z^{- 1/k + 1} \frac{\partial}{\partial x}
\Delta_k^x (z)^{-1}  &=& k z^{- 1/k + 1} 
\frac{\partial}{\partial z} \Delta_k^x (z)^{-1} . 
\label{second identity in rep}
\end{eqnarray}
\end{prop}

\begin{proof}  In $\C[x, x^{-1}] [[z^{1/k}, z^{-1/k}]]$, we have
\begin{eqnarray*}
& & \hspace{-.5in} 
- \Delta_k^x (z) \frac{\partial}{\partial x} \cdot x + \frac{1}{k} 
z^{1/k - 1} \frac{\partial}{\partial x} 
\Delta_k^x (z) \cdot x =  \\
&=& \! -1 + \frac{1}{k} z^{1/k - 1} \frac{\partial}{\partial x}
\Biggl( k z^{- 1/k + 1} z^{- (1/k) \Lo} z^{1/k}
\exp \Biggl( - \sum_{j \in \ZZ} a_j \Lx \Biggr) \cdot x \Biggr)\\ 
&=& \! -1 + \frac{1}{k} z^{1/k - 1} \frac{\partial}{\partial x}
\left( k z^{- 1/k + 1} z^{- (1/k) \Lo} f(x) \right)\\
&=& \! -1 + \frac{1}{k} z^{1/k - 1} k z^{- 1/k + 1}
\frac{\partial}{\partial x} f(z^{- 1/k} x)\\ 
&=& \! -1 + z^{- 1/k} f'(z^{- 1/k} x)\\ 
&=& \! -1 + (1+ z^{- 1/k} x)^{k-1}\\ 
&=& \! -1 + (1+ z^{- 1/k} x)^k - z^{- 1/k} x (1+ z^{- 
1/k} x)^{k-1}\\
&=& \! (-1 + k) z^{- 1/k} \Biggl( \frac{ z^{1/k}}{k}
(1 + z^{- 1/k}x)^k - \frac{ z^{1/k}}{k} \Biggr) \\
& & \hspace{1.7in} + \;  k z^{- 1/k + 1} 
\frac{\partial}{\partial z} \Biggl( \frac{ z^{1/k}}{k}(1 + 
z^{- 1/k}x)^k - \frac{ z^{1/k}}{k} \Biggr)\\
&=& \! k( - \frac{1}{k} + 1) z^{-1/k} f(z^{-1/k}x) + k
z^{-1/k + 1} \frac{\partial}{\partial z} f(z^{-1/k}x)
\\
&=& \! \frac{\partial}{\partial z} \left( k z^{-1/k + 1}
f(z^{-1/k}x) \right)\\
&=& \! \frac{\partial}{\partial z} \left( k z^{-1/k + 1}
z^{-(1/k) \Lo} f(x) \right) \\
&=& \! \frac{\partial}{\partial z} \Biggl( k z^{-1/k + 1}
z^{-(1/k) \Lo} z^{1/k} \exp \Biggl(- \sum_{j \in \ZZ} a_j
\Lx \Biggr) \cdot x \Biggr) \\
&=& \! \frac{\partial}{\partial z} \Biggl( k z^{-1/k + 1} \exp 
\Biggl(- \sum_{j \in \ZZ} a_j z^{- j/k} \Lx \Biggr) \cdot x \Biggr)\\
&=& \! \frac{\partial}{\partial z} \Delta_k^x (z) \cdot x .
\end{eqnarray*}

Since $\Delta_k^x (z) \cdot  x^n$ is well defined in $\C
[x, x^{-1}] [[ z^{1/k}, z^{- 1/k}]]$ for all 
$n \in \Z$, by Proposition 2.1.7 in \cite{H}, we have
\[\Delta_k^x (z_2) \cdot  x^n = (\Delta_k^x (z_2) \cdot  x)^n \]
for all $n \in \Z$.  Therefore
\begin{eqnarray*}
& & \hspace{-.7in} 
- \Delta_k^x (z) \frac{\partial}{\partial x} \cdot x^n  + 
\frac{1}{k} z^{1/k - 1} \frac{\partial}{\partial x} 
\Delta_k^x (z) \cdot x^n = \\
&=& \! -n  \Delta_k^x (z) \cdot x^{n - 1}
\frac{\partial}{\partial x} \cdot x  + \frac{1}{k} z^{1/k - 1} 
\frac{\partial}{\partial x} \left(\Delta_k^x (z)  \cdot x \right)^n \\
&=& \! -n \left(\Delta_k^x (z) \cdot x^{n - 1} \right) \left(\Delta_k^x (z)
\frac{\partial}{\partial x} \cdot x \right) \\
& & \hspace{1.6in} + \; n \frac{1}{k} z^{1/k - 1} 
\left( \Delta_k^x (z) \cdot x \right)^{n - 1} \frac{\partial}{\partial x} 
\left( \Delta_k^x (z) \cdot x  \right) \\
\hspace{.3in} &=& n \left( \Delta_k^x (z) \cdot x \right)^{n - 1} \left(- 
\Delta_k^x (z)  
\frac{\partial}{\partial x} \cdot x  + \frac{1}{k} z^{1/k - 1} 
\frac{\partial}{\partial x} \left( \Delta_k^x (z) \cdot x  \right) \right)\\
&=& \! n \left( \Delta_k^x (z)  \cdot x 
\right)^{n - 1} \left( \frac{\partial}{\partial z} \left( \Delta_k^x (z) 
 \cdot x \right) \right)\\
&=& \! \frac{\partial}{\partial z} \left( \Delta_k^x (z) \cdot x \right)^n \\
&=& \! \frac{\partial}{\partial z} \left( \Delta_k^x (z) \cdot x^n \right)\\
\end{eqnarray*}
for all 
$n \in \Z$.  Equation (\ref{identity in rep}) follows by linearity.

Similarly, in $\C[x, x^{-1}] [[z^{1/k}, z^{-1/k}]]$, we have
\begin{eqnarray*}
& & \hspace{-.4in} 
- \Delta_k^x (z)^{-1} \frac{\partial}{\partial x} \cdot x + k 
z^{- 1/k + 1} \frac{\partial}{\partial x} 
\Delta_k^x (z)^{-1} \cdot x = \\
&=& \! -1 + k z^{- 1/k + 1} \frac{\partial}{\partial x}
\Biggl( z^{(1/k - 1) \Lo} k^{- \Lo} \exp \Biggl( 
\sum_{j \in \ZZ} a_j z^{- j/k} \Lx \Biggr) \cdot x \Biggr)\\
&=& \! -1 + k z^{- 1/k + 1} \frac{\partial}{\partial x}
\left( z^{(1/k - 1) \Lo} k^{- \Lo}
z^{1/k} f^{-1}(x) \right)\\
&=& \! -1 + k z \frac{\partial}{\partial x}
\left( f^{-1}( z^{1/k - 1 } k^{-1} x) \right)\\
&=& \! -1 + k z \frac{\partial}{\partial x}
\left( (1 + z^{-1} x )^{1/k} - 1\right)\\
&=& \! -1 + (1 + z^{-1} x )^{1/k - 1}\\
&=& \! (1 + z^{-1} x )^{1/k} - 1 - z^{- 1}x 
(1 + z^{-1} x )^{1/k - 1} \\
&=& \! k z^{- 1/k + 1} \left( \! \Bigl( \frac{\partial}{\partial z} 
z^{1/k} \Bigr) \! \left( (1 + z^{-1} x )^{1/k} 
- 1 \right) + z^{1/k} \frac{\partial}{\partial z} \! 
\left( (1 + z^{-1} x )^{1/k} - 1\right) \! \right) \\
&=& \! k z^{- 1/k + 1} \frac{\partial}{\partial z} \left(   
z^{1/k} f^{-1} (z^{1/k - 1}k^{-1}  x) \right)\\
&=& \! k z^{- 1/k + 1} \frac{\partial}{\partial z} 
\left( z^{(1/k - 1) \Lo} k^{- \Lo} z^{1/k} 
f^{-1} (x) \right)\\
&=& \! k z^{- 1/k + 1} \frac{\partial}{\partial z} 
\Biggl( z^{(1/k - 1) \Lo} k^{- \Lo} \exp \Biggl(
\sum_{j \in \ZZ} a_j z^{- j/k} \Lx \Biggr) 
\cdot x \Biggr)\\
&=& \! k z^{- 1/k + 1} \frac{\partial}{\partial z} 
\Delta_k^x (z)^{-1} \cdot x .
\end{eqnarray*}
The proof of identity (\ref{second identity in rep}) on $x^n$ for 
$n \in \Z$ is analogous to the proof of identity 
(\ref{identity in rep}) on $x^n$ for $n \in \Z$.  Identity 
(\ref{second identity in rep}) then follows by linearity.  
\end{proof}

Let $\mathfrak{L}$ be the Virasoro algebra with basis $L_j$, $j \in \Z$,
and central charge $d \in \C$.  The above identity can be thought of as 
an identity for the representation of the Virasoro algebra on 
$\C[x, x^{-1}]$ given by $L_j \mapsto -\Lx$, for $j \in \Z,$ with 
central charge  equal to zero.  We want to prove the corresponding 
identity for certain other representations of the Virasoro algebra, in 
particular for vertex operator algebras.  We do this by following the 
method of proof used in Chapter 4 of \cite{H}.  Letting $\kappa$ be
another formal variable commuting with $z$ and $\mathfrak{L}$, we first 
prove the identity in $\mathcal{ U}_\Pi (\mathfrak{L})  [[z^{1/k}, 
z^{-1/k}]][[\kappa, \kappa^{-1}]]$ where $\mathcal{ U}_\Pi 
(\mathfrak{L})$ is a certain extension of the universal enveloping 
algebra for the Virasoro algebra, and then letting $\kappa = k$, the 
identity will follow in $(\End \; V) [[z^{1/k}, 
z^{-1/k}]]$ where $V$ is a certain type of module for the 
Virasoro algebra.

We want to construct an extension of $\mathcal{ U}(\mathfrak{L})$, the 
universal enveloping algebra for the Virasoro algebra, in which 
$\kappa^{-L_0}$ and $z^{(1/k - 1)L_0}$ can be defined.  Let 
$V_\Pi$ be a vector space over $\mathbb{C}$ with basis $\{P_j : j \in 
\mathbb{Z}\}$. Let $\mathcal{ T}(\mathfrak{L} \oplus V_\Pi)$ be the tensor 
algebra generated by the direct sum of $\mathfrak{L}$ and $V_\Pi$, and 
let $\mathcal{ I}$ be the ideal of $\mathcal{ T}(\mathfrak{L} \oplus V_\Pi)$ 
generated by
\begin{eqnarray*}
\lefteqn{\Bigl\{L_i \otimes L_j - L_j \otimes L_i - [L_i,L_j], \; L_i \otimes d - d 
\otimes L_i, \; P_i \otimes P_j - \delta_{ij} P_i,  \Bigr. }   \\
& & \hspace{1.5in} P_i \otimes d -
d \otimes P_i, \; \Bigl. P_i \otimes L_j - L_j \otimes P_{i + j} : \; i,j \in 
\mathbb{Z}\Bigr\}.
\end{eqnarray*}
Define $\mathcal{ U}_\Pi (\mathfrak{L}) = \mathcal{ T}(\mathfrak{L} \oplus V_\Pi)/
\mathcal{ I}$.  For any formal variable $z$ and for $n\in\mathbb{Z}$, we define
\[z^{nL_0} = \sum_{j \in \mathbb{Z}} P_j z^{nj} \in \mathcal{ U}_\Pi(\mathfrak{L})
[[z,z^{-1}]] .\]
Note $\kappa^{-L_0}$ and $z^{(1/k - 1)L_0}$ are well-defined elements
of $U_\Pi(\mathfrak{L})[[z^{1/k},z^{-1/k}]][[\kappa, \kappa^{-1}]]$.

In $\mathcal{ U}_\Pi (\mathfrak{L})[[z^{1/k}, z^{-1/k}]]
[[\kappa, \kappa^{-1}]]$, define 
\[\Delta_k^\mathfrak{L} (z) = \exp \Biggl( \sum_{j \in \ZZ} 
a_j z^{- j/k} L_j \Biggr) \kappa^{-L_0} z^{\left( 1/k - 1 \right) 
L_0} .\]

\begin{prop}\label{in universal}
In $\mathcal{ U}_\Pi (\mathfrak{L})[[z^{1/k}, z^{-1/k}]]
[[\kappa, \kappa^{-1}]]$, we have
\begin{eqnarray}
\Delta_k^\mathfrak{L} (z) L_{-1} - \frac{1}{\kappa} z^{1/k - 1} L_{-1} 
\Delta_k^\mathfrak{L} (z) 
&= &\frac{\partial}{\partial z}\Delta_k^\mathfrak{L} (z) ,
\label{identity in universal}\\
\Delta_k^\mathfrak{L} (z)^{-1} L_{-1} - \kappa z^{-1/k + 1} L_{-1} 
\Delta_k^\mathfrak{L} (z)^{-1} 
&=& k z^{-1/k + 1} \frac{\partial}{\partial z}
\Delta_k^\mathfrak{L} (z)^{-1} . 
\label{second identity in universal}
\end{eqnarray} 
\end{prop}

\begin{proof}  In $\mathcal{ U}_\Pi (\mathfrak{L})[[z^{1/k}, 
z^{-1/k}]][[\kappa, \kappa^{-1}]]$, we have
\begin{eqnarray*}
& & \hspace{-.4in}
\Delta_k^\mathfrak{L} (z) L_{-1} - \frac{1}{\kappa} z^{1/k - 1} L_{-1} 
\Delta_k^\mathfrak{L} (z) = \\
&=& \! \frac{1}{\kappa} z^{1/k - 1} \left[ e^{\sum_{j \in \ZZ} 
a_j z^{- j/k} L_j}, L_{-1} \right] \kappa^{- L_0} z^{\left(1/k - 1 
\right) L_0} \\
&=& \! \frac{1}{\kappa} z^{1/k - 1} \sum_{n \in \ZZ} \frac{1}{n!}\Biggl( 
\sum_{j_1,...,j_n \in \ZZ} a_{j_1} \cdots a_{j_n}  z^{-(j_1 + \cdots + j_n)/k} 
\Biggr. \\ 
& & \hspace{.6in} \Biggl. \biggl( \sum_{i = 1,...,n} L_{j_1} L_{j_2} \cdots
L_{j_{i-1}}[L_{j_i},L_{-1}] L_{j_{i + 1}} \cdots L_{j_n} \biggr)\!\Biggr)
\kappa^{- L_0} z^{\left(1/k - 1 \right) L_0} 
\end{eqnarray*} 
which is a well-defined element of $\mathcal{ U}_\Pi (\mathfrak{L})
[[z^{1/k}, z^{-1/k}]][[\kappa, \kappa^{-1}]]$ involving 
only elements $L_j$ with $j \in {\N}$.  The right-hand side of 
(\ref{identity in universal}) also involves only $L_j$ for $j \in {\N}$.  
Thus comparing with the identity (\ref{identity in rep}) for the 
representation $L_j \mapsto - \Lx$, the identity (\ref{identity in 
universal}) must hold.  The proof of (\ref{second identity in 
universal}) is analogous.  \end{proof}

Let $V$ be a module for the Virasoro algebra satisfying $V = \coprod_{n 
\in \mathbb{Z}} V_n$.  For $j \in \mathbb{Z}$, let $L(j) \in \End \; V$  
and $c \in \C$ be the representation images of $L_j$ and $d$, 
respectively, for the Virasoro algebra.  Assume that for $v \in V_n$, 
we have $L(0)v = nv$.  For any formal variable $z$, define $z^{jL(0)} 
\in (\mathrm{End} \; V)[[z,z^{-1}]]$ by
\[z^{jL(0)}v = z^{jn}v \]
for $v \in V_n$.  Or equivalently, let $P(n) : V \rightarrow V_n$
be the projection {}from $V$ to the homogeneous subspace of weight $n$ for
$n \in \mathbb{Z}$.  Then
\[z^{jL(0)}v = \sum_{n \in \mathbb{Z}} z^{jn} P(n) v \]
for $v \in V$.  The elements $P(n) \in \mathrm{End} \; V$ can be thought
of as the representation images of $P_n$ in the algebra $\mathcal{ U}_\Pi
(\mathfrak{L})$.

Note that for $k$ a positive integer, $k^{-L(0)}$ is a well-defined element of 
$\mathrm{End} \; V$ and $z^{(1/k - 1)L(0)}$ is a well-defined element
of $(\mathrm{End} \; V)[[z^{1/k},z^{-1/k}]]$.

In $(\End \; V)[[z^{1/k}, z^{-1/k}]]$, define 
\begin{equation}\label{Delta for a module}
\Delta_k^V (z) = \exp \Biggl( \sum_{j \in \ZZ} a_j z^{- j/k} L(j) 
\Biggr) k^{-L(0)} z^{\left( 1/k - 1 \right) L(0)} .
\end{equation}
{}From Proposition 4.1.1 in \cite{H} and Proposition \ref{in universal}, we 
obtain the following corollary.

\begin{cor}\label{c2.5}
In $(\End \; V)[[z^{1/k}, z^{-1/k}]]$, we have
\begin{eqnarray}
\Delta_k^V (z) L(-1) - \frac{1}{k} z^{1/k - 1} L(-1) 
\Delta_k^V (z) \!
&=& \! \frac{\partial}{\partial z}\Delta_k^V (z) , \label{identity in voa}\\ 
\Delta_k^V (z)^{-1} L(-1) - k z^{- 1/k + 1} L(-1) 
\Delta_k^V (z)^{-1}  \! 
&=& \! k z^{- 1/k + 1} \frac{\partial}{\partial z}
\Delta_k^V (z)^{-1}  . \label{second identity in voa}
\end{eqnarray}
In particular, the identities hold for $V$ being any vertex operator
algebra. 

\end{cor}

\section{The twisted sector for $g = (1 2 \cdots k)$}
\setcounter{equation}{0}

We first review the definitions of weak, weak admissible and
ordinary $g$-twisted modules for a vertex operator algebra $V$
and an automorphism $g$ of $V$ of finite order $k$ (cf. \cite{DLM1}--
\cite{DLM3}). 

Let $(V, Y, \mathbf{1}, \omega)$ be a vertex operator algebra.  A {\em
weak} $g$-twisted $V=(V,Y_M)$-module is a $\C$-linear space $M$
equipped with a linear map $V\to (\End M)[[z^{1/k},z^{-1/k}]]$, given
by $v\mapsto Y_M(v,z)=\sum_{n\in\Q}v_nz^{-n-1}$, such that for $u,v\in
V$ and $w\in M$ the following hold: (1) $v_mw=0$ if $m$ is
sufficiently large; (2) $Y_M({\bf 1},z)=1;$ (3) $Y_M(v,z)=\sum_{n\in
r/k+\Z}v_nz^{-n-1}$ for $v\in V^r$ where $V^r=\{v\in V|gv=e^{-2\pi
ir/k}v\};$ (4) the twisted Jacobi identity holds: for $u\in V^r$
\[z^{-1}_0\delta\left(\frac{z_1-z_2}{z_0}\right)
Y_M(u,z_1)Y_M(v,z_2)-z^{-1}_0\delta\left(\frac{z_2-z_1}{-z_0}\right)
Y_M(v,z_2)Y_M(u,z_1) \]
\begin{equation}\label{g3.13}
= z_2^{-1}\left(\frac{z_1-z_0}{z_2}\right)^{-r/k}
\delta\left(\frac{z_1-z_0}{z_2}\right) Y_M(Y(u,z_0)v,z_2). 
\end{equation}
It can be shown (cf. Lemma 2.2 of \cite{DLM1}, and \cite{DLM2}) that 
$Y_M(\o,z)$ has component operators which satisfy the Virasoro algebra 
relations and $Y_M(L(-1)u,z)=\frac{d}{dz}Y_M(u,z)$. If we take $g=1$, 
then we obtain a weak $V$-module.

A {\em weak admissible} $g$-twisted $V$-module is a weak $g$-twisted
$V$-module $M$ which carries a $\frac{1}{k}{\Z}_{+}$-grading
\begin{equation}\label{m3.12}
M=\oplus_{n\in\frac{1}{k}\Z_{+}}M(n)
\end{equation}
such that $v_mM(n)\subseteq M(n+\wt \; v-m-1)$ for homogeneous $v\in V.$
We may assume that $M(0)\ne 0$ if $M\ne 0$.  If $g=1,$ we have a
weak admissible $V$-module.

\begin{rem}
Above we used the term ``weak admissible $g$-twisted module''
whereas in most of the literature (cf. \cite{DLM1}, \cite{Z}) the term
``admissible $g$-twisted module'' is used for this notion.  We used
the qualifier ``weak'' to stress that these are indeed only weak
modules and in general are not ordinary modules.  However, for the
sake of brevity, we will now drop the qualifier ``weak''.
\end{rem}

An (ordinary) $g$-twisted $V$-module is a weak $g$-twisted $V$-module $M$ 
graded by $\C$ induced by the spectrum of $L(0).$ That is, we have
\begin{equation}\label{g3.14}
M=\coprod_{\lambda \in{\C}}M_{\lambda} 
\end{equation}
where $M_{\l}=\{w\in M|L(0)w=\l w\}.$ Moreover we require that $\dim
M_{\l}$ is finite and $M_{n/k +\l}=0$ for fixed $\l$ and for all
sufficiently small integers $n.$ If $g=1$ we have an ordinary
$V$-module.

The vertex operator algebra $V$ is called $g$-{\em rational} if every
admissible $g$-twisted $V$-module is completely reducible, i.e., a
direct sum of irreducible admissible $g$-twisted modules. It was
proved in \cite{DLM2} that if $V$ is $g$-rational then: (1) every
irreducible admissible $g$-twisted $V$-module is an ordinary
$g$-twisted $V$-module; and (2) $V$ has only finitely many isomorphism
classes of irreducible admissible $g$-twisted modules.

Now we turn our attention to tensor product vertex operator algebras.
Let $V=(V,Y,{\bf 1},\omega)$ be a vertex operator algebra and $k$ a
fixed positive integer as in Section 2. Then $V^{\otimes k}$ is also a
vertex operator algebra (see \cite{FHL}), and the permutation group $S_k$
acts naturally on $\vtk$ as automorphisms. Let $g=(1 2\cdots k)$.  In
this section we construct a functor $T_g^k$ {}from the category of weak
$V$-modules to the category of weak $g$-twisted modules for $\vtk.$ We
do this by first defining $g$-twisted vertex operators on a weak
$V$-module $M$ for a set of generators which are mutually local (see
\cite{Li2}).  These $g$-twisted vertex operators generate a local system 
which is a vertex algebra.  We then construct a homomorphism of vertex 
algebras {}from $V^{\otimes k}$ to this local system which thus gives a
weak $g$-twisted $V^{\otimes k}$-module structure on $M$.

For $v\in V$ denote by $v^j\in\vtk$ the vector whose $j$-th tensor
factor is $v$ and whose other tensor factors are ${\bf 1}$.  Then
$gv^j=v^{j+1}$ for $j=1,...,k$ where $k+1$ is understood to be 1. Let
$W$ be a weak $g$-twisted $\vtk$-module, and let $\eta=e^{-2\pi i/k}$.
Then it follows immediately {}from the definition of twisted module that
the $g$-twisted vertex operators on $W$ satisfy
$$Y_g(v^{j+1},z)=\lim_{z^{1/k}\to \eta^{-j}z^{1/k}}Y_g(v^1,z).$$
Since $\vtk$ is generated by $v^j$ for $v\in V$ and $j=1,...,k,$ 
the vertex operators $Y_g(v^1,z)$ for $v\in V$ determine all
the vertex operators $Y_g(u,z)$ on $W$ for any $u\in\vtk.$ This
observation is very important in our construction of twisted
sectors.

Let $u,v\in V$.  Then by (\ref{g3.13}) the twisted Jacobi identity for
$Y_g(u^1,z_1)$ and $Y_g(v^1,z_2)$ is
\[z^{-1}_0\delta\left(\frac{z_1-z_2}{z_0}\right)
Y_g(u^1,z_1)Y_g(v^1,z_2)-z^{-1}_0\delta\left(\frac{z_2-z_1}{-z_0}\right)
Y_g(v^1,z_2)Y_g(u^1,z_1)\]
\begin{equation}\label{k1}
=\frac{1}{k}z_2^{-1}\sum_{j=0}^{k-1}\delta\Biggl(\eta^j\frac{(z_1-z_0)^
{1/k}}{z_2^{1/k}}\Biggr)Y_g(Y(g^ju^1,z_0)v^1,z_2)
\end{equation} 
(cf. \cite{Le2}, \cite{D}).  Since $g^ju^1=u^{j+1}$, we see that
$Y(g^ju^1,z_0)v^1$ only involves nonnegative integer powers of $z_0$
unless $j=0\ (\mod \; k).$ 
Thus
\begin{equation}\label{k2}
[Y_g(u^1,z_1),Y_g(v^1,z_2)] \; = \; \Res_{z_0}\frac{1}{k}z_2^{-1}\delta
\Biggl(\frac{(z_1-z_0)^{1/k}}{z_2^{1/k}}\Biggr)Y_g(Y(u^1,z_0)v^1,z_2).
\end{equation}
This shows that the component operators of $Y_g(u^1,z)$ for $u\in V$ on
$W$ form a Lie algebra. 

Now let $M=(M,Y)$ be a weak $V$-module. For $u\in V$ and $\D(z) = 
\D^V(z)$ given by (\ref{Delta for a module}), define
\[\bar Y(u,z)=Y(\D(z)u,z^{1/k}) .\] 
When we put a weak $g$-twisted $\vtk$-module structure on $M$ this
$\bar{Y}(u,z)$ will be the twisted vertex operator acting on $M$
associated to $u^1$.  Here we give several examples of $\bar Y(u,z)$.

If $u\in V_n$ is a highest weight vector then
$\Delta_k(z)u=k^{-n}z^{(1/k-1)n}u$ and
\[\bY(u,z)=k^{-n}z^{(1/k-1)n}Y(u,z^{1/k}).\]
In particular, if $n=1$ we have 
\[\bY(u,z)=k^{-1}z^{1/k-1}Y(u,z^{1/k}).\]
This case is important in the study of symmetric orbifold theory 
for the vertex operator algebras associated to affine Lie algebras.
Now we take $u=\omega.$ Recall that $a_2= (k^2-1)/12$.  Thus
\begin{eqnarray*}
\Delta_k(z)\omega &=& \frac{z^{2(1/k - 1)}}{k^2}\Bigl(\omega + a_2
\frac{c}{2}z^{-2/k} \Bigr)\\
&=& \frac{z^{2(1/k - 1)}}{k^2}\Bigl(\omega + \frac{(k^2-1)c}{24}
z^{-2/k} \Bigr)
\end{eqnarray*}
where $c$ is the central charge. Therefore 
\begin{equation}\label{sun1}
\bY(\omega,z)=\frac{z^{2(1/k - 1)}}{k^2}Y(\omega,z^{1/k})
+\frac{(k^2-1)c}{24k^2}z^{-2}.
\end{equation}

We next study the properties of the operators $\bY(u,z).$

\begin{lem}\label{l3.1} For $u\in V$
\[\bY(L(-1)u,z)=\frac{d}{dz}\bY(u,z).\]
\end{lem}

\begin{proof} By Corollary \ref{c2.5}, we have
\begin{eqnarray*}
\bY(L(-1)u,z) &=& Y(\D(z)L(-1)u,z^{1/k})\\
&=& Y(\frac{d}{dz}\D(z)u,z^{1/k}) + \frac{1}{k}z^{1/k - 1}Y(L(-1)\D(z)u,z^{1/k})\\
&=& Y(\frac{d}{dz}\D(z)u,z^{1/k}) + \left.\frac{1}{k}z^{1/k - 1}\frac{d}{dx}
Y(\D(z)u,x) \right|_{x=z^{1/k}}\\
&=& Y(\frac{d}{dz}\D(z)u,z^{1/k}) + \left.\frac{d}{dx}Y(\D(z)u,x^{1/k})
\right|_{x=z}\\
&=& \frac{d}{dz}Y(\D(z)u,z^{1/k})\\
&=& \frac{d}{dz}\bY(u,z)
\end{eqnarray*}
as desired. \end{proof}

In the proof of the following lemma and again later on, we will need
some properties of the $\delta$-function.  We first note that {}from
Proposition 8.8.22 of \cite{FLM3}, for $p \in \mathbb{Z}$ we have
\begin{equation}\label{delta function1}
z_2^{-1}\left(\frac{z_1-z_0} {z_2}\right)^{-p/k}
\delta\left(\frac{z_1-z_0} {z_2}\right) =
z_1^{-1}\left(\frac{z_2+z_0} {z_1}\right)^{p/k}\delta
\left(\frac{z_2+z_0} {z_1}\right) ,
\end{equation}
and it is easy to see that  
\begin{equation}\label{delta function2}
\sum_{p=0}^{k-1}\left(\frac{z_1-z_{0}}{z_2}\right)^{p/k}
z_2^{-1}\delta\left(\frac{z_1-z_0}{z_2}\right)=
z_2^{-1}\delta\Biggl(\frac{(z_1-z_0)^{1/k}}{z_2^{1/k}}\Biggr).
\end{equation}
Therefore, we have the $\delta$-function identity
\begin{equation}\label{delta function3}
z_2^{-1} \delta \Biggl( \frac{(z_1 - z_0)^{1/k}}{z_2^{1/k}} \Biggr) = 
z_1^{-1} \delta \Biggl( \frac{(z_2 + z_0)^{1/k}}{z_1^{1/k}} \Biggr).
\end{equation}

\begin{lem}\label{l3.2} For $u,v\in V$ 
\[ [\bY(u,z_1),\bY(v,z_2)] \; = \; \Res_{z_0}\frac{1}{k}z_2^{-1}
\delta\Biggl(\frac{(z_1-z_0)^{1/k}}{z_2^{1/k}}\Biggr)\bY(Y(u,z_0)v,z_2).\]
\end{lem}

\begin{proof} Replacing $Y(u,z_1)$ and $Y(v,z_2)$ by $Y(\D(z_1)u,z_1^{1/k})$ and 
$Y(\D(z_2)v,z_2^{1/k})$, respectively, in the commutator
formula
\begin{equation}\label{commutator for Lemma 3.3}
[Y(u,z_1),Y(v,z_2)] \; = \; \Res_{x}z_2^{-1}\delta\left(\frac{z_1-x}{z_2}
\right)Y(Y(u,x)v,z_2) 
\end{equation}
which is a consequence of the Jacobi identity on $M,$ we have
\begin{equation}\label{substitution equation}
[\bY(u,z_1),\bY(v,z_2)] \; = \; \Res_{x}z_2^{-1/k}\delta\Biggl(
\frac{z_1^{1/k}-x}{z_2^{1/k}}\Biggr)
Y(Y(\D(z_1)u,x)\D(z_2)v,z_2^{1/k}).
\end{equation}

We want to make the change of variable $x = z_1^{1/k}-(z_1-z_0)^{1/k}$
where by $z_1^{1/k}-(z_1-z_0)^{1/k}$ we mean the power series expansion
in positive powers of $z_0$.  In this case, we note that for $n \in \mathbb{Z}$
\begin{eqnarray*}
\hspace{.8in} & & \hspace{-1.2in} \left. (z_1^{1/k} - x)^n \right|_{x = z_1^{1/k}-(z_1-z_0)^{1/k}} =
\\ 
&=& \!  \sum_{m\in \mathbb{N}} \binom{n}{m} (-1)^m z_1^{n/k - m/k} 
\biggl( - \! \sum_{l \in \mathbb{Z}_+} \binom{1/k}{l} z_1^{1/k - l} (-1)^l z_0^l \biggr)^m\\
&=& \!  \sum_{m\in \mathbb{N}} \binom{n}{m} (-1)^m z_1^{n/k} 
\biggl( - \sum_{l \in \mathbb{Z}_+} \binom{1/k}{l} \Bigl(\frac{-z_0}{z_1}\Bigr)^l \biggr)^m\\
&=& \!  z_1^{n/k} \biggl( 1 + \sum_{l \in \mathbb{Z}_+} \binom{1/k}{l}
\Bigl(\frac{-z_0}{z_1}\Bigr)^l \biggr)^n\\ &=& \! \! z_1^{n/k}\Bigl(1 -
\frac{z_0}{z_1}\Bigr)^{n/k}\\ &=& \! \! (z_1 - z_0)^{n/k}.
\end{eqnarray*}

Thus substituting $x =  z_1^{1/k}-(z_1-z_0)^{1/k}$ into 
\[z_2^{-1/k}\delta\Biggl(\frac{z_1^{1/k}-x}{z_2^{1/k}}\Biggr)
Y(Y(\D(z_1)u,x)\D(z_2)v,z_2^{1/k}) \]
we have a well-defined power series given by
\[\delta\Biggl(\frac{(z_1-z_0)^{1/k}}{z_2^{1/k}}\Biggr)
 Y(Y(\D(z_1)u, z_1^{1/k}-(z_1-z_0)^{1/k})\D(z_2)v,z_2^{1/k}) .\] 

Let $f(z_1,z_2,x)$ be a complex analytic function in $z_1, z_2$, and
$x$, and let $h(z_1,z_2,z_0)$ be a complex analytic function in $z_1,
z_2$, and $z_0$. Then if $f(z_1,z_2,h(z_1,z_2,z_0))$ is well defined,
and thinking of $z_1$ and $z_2$ as fixed, i.e., considering 
$f(z_1,z_2,h(z_1,z_2,z_0))$ as a Laurent series in $z_0$, by the
residue theorem of complex analysis, we have
\begin{equation}\label{residue change of variables}
\Res_x f(z_1,z_2,x)=\Res_{z_0} \left( \frac{\partial}{\partial z_0}
h(z_1,z_2,z_0) \right)f(z_1,z_2,h(z_1,z_2,z_0))
\end{equation} 
which of course remains true for $f$ and $h$ formal power series in
their respective variables.  Thus making the change of variable
$x= h(z_1,z_2,z_0) = z_1^{1/k}-(z_1-z_0)^{1/k}$, using (\ref{substitution
equation}), (\ref{residue change of variables}), the 
$\delta$-function identity (\ref{delta function3}) and Proposition 
\ref{psun1}, we obtain
\begin{eqnarray*}
& & \hspace{-.35in} [\bY(u,z_1),\bY(v,z_2)] = \\  
&=& \! \Res_{z_0}\frac{1}{k}z_2^{-1/k} (z_1-z_0)^{1/k-1}\delta\Biggl(
\frac{(z_1-z_0)^{1/k}}{z_2^{1/k}}\Biggr) \\
& & \! \hspace{1.7in} Y(Y(\D(z_1)u,z_1^{1/k}-(z_1-z_0)^{1/k})
\D(z_2)v,z_2^{1/k})\\
&=& \! \Res_{z_0}\frac{1}{k} z_2^{-1} \delta\Biggl(
\frac{(z_1-z_0)^{1/k}}{z_2^{1/k}}\Biggr)  
Y(Y(\D(z_1)u,z_1^{1/k}-(z_1-z_0)^{1/k})\D(z_2)v,z_2^{1/k})\\
&=& \! \Res_{z_0}\frac{1}{k}z_1^{-1} \delta\Biggl(
\frac{(z_2+z_0)^{1/k}}{z_1^{1/k}}\Biggr)  
Y(Y(\D(z_1)u,z_1^{1/k}-(z_1-z_0)^{1/k})\D(z_2)v,z_2^{1/k})\\
&=& \! \Res_{z_0}\frac{1}{k}z_2^{-1}\delta\Biggl(
\frac{(z_1-z_0)^{1/k}}{z_2^{1/k}}\Biggr)  \\
& &\! \hspace{1.4in}  Y(Y(\D(z_2+z_0)u,(z_2+z_0)^{1/k}-z_2^{1/k})
\D(z_2)v,z_2^{1/k})\\
&=& \! \Res_{z_0}\frac{1}{k}z_2^{-1}\delta\Biggl(\frac{(z_1-z_0)^
{1/k}}{z_2^{1/k}}\Biggr)Y(\D(z_2)Y(u,z_0)v,z_2)\\
&=& \! \Res_{z_0}\frac{1}{k}z_2^{-1}\delta\Biggl(\frac{(z_1-z_0)^
{1/k}}{z_2^{1/k}}\Biggr)\bY(Y(u,z_0)v,z_2),
\end{eqnarray*}
as desired. \end{proof}

We are now in a position to put a weak $g$-twisted $\vtk$-module
structure on $M.$ For $u\in V$ set
\begin{equation}\label{define g-twist}
Y_g(u^1,z)=\bY(u,z)  \quad \mbox{and} \quad   
Y_g(u^{j+1},z)=\lim_{z^{1/k}\to \eta^{-j} z^{1/k}}Y_g(u^1,z).
\end{equation}
Note that $Y_g(u^j,z)=\sum_{p=0}^{k-1}Y_g^p(u^j,z)$
where $Y_g^p(u^j,z)=\sum_{n\in p/k+\Z}u^j_nz^{-n-1}$.

\begin{lem}\label{l3.3} Let $u,v\in V.$ Then
\begin{equation}\label{3.1}
[Y_g(u^i,z_1),Y_g(v^j,z_2)] \; = \; \Res_{z_0}\frac{1}{k}z_2^{-1}\delta
\Biggl(\frac{\eta^{j-i}(z_1-z_0)^{1/k}}{z_2^{1/k}}\Biggr)Y_g((Y(u,z_0)v)^j,z_2)
\end{equation}
where $(Y(u,z_0)v)^j=\sum_{n\in\Z}(u_nv)^jz_0^{-n-1},$
and
\begin{multline}\label{3.2}
[Y_g^p(u^i,z_1),Y_g(v^j,z_2)] \\
= \Res_{z_0}\frac{1}{k}z_2^{-1}\eta^{(j-i)p}\left(\frac{z_1-z_0}
{z_2}\right)^{-p/k}\delta\left(\frac{z_1-z_0}
{z_2}\right)Y_g((Y(u,z_0)v)^j,z_2).
\end{multline}
\end{lem}

\begin{proof} By Lemma \ref{l3.2}, equation (\ref{3.1}) holds if $i=j=1.$
Replacing $z_1^{1/k}$ and $z_2^{1/k}$ by $\eta^{-i+1}z_1^{1/k}$ and
$\eta^{-j+1}z_2^{1/k}$, respectively, we obtain equation (\ref{3.1})
for any $i,j = 1,...,k.$ Equation (\ref{3.2}) is a direct consequence of
(\ref{3.1}). \end{proof}

By Lemma \ref{l3.3} for $u,v\in V$, there exists a positive integer
$N$ such that 
\begin{equation}\label{a3.3}
[Y_g(u^i,z_1), Y_g(v^j,z_2)](z_1-z_2)^N=0.
\end{equation}
Letting $z^{1/k}$ go to $\eta^{-i+1}z^{1/k}$ in Lemma \ref{l3.1}, we have
$$Y_g(L(-1)u^i,z)=\frac{d}{dz}Y_g(u^i,z).$$ Thus the operators
$Y_g(u^i,z)$ for $u \in V$, and $i=1,...,k$ are mutually local and
generate a local system $A$ in the sense of \cite{Li2}.  Let $\sigma$ be a
map {}from $A$ to $A$ such that $\sigma Y_g(u^i,z)=Y_g(u^{i+1},z)$ for
$u\in V$ and $i=1,...,k$.  By Theorem 3.14 of \cite{Li2}\footnote{There is a
typo in the statement of Theorem 3.14 in \cite{Li2}.  The $V$ in the theorem
should be $A$.  That is, the main result of the theorem is that the
local system $A$ of the theorem has the structure of a vertex
superalgebra.}, the local system $A$ generates a vertex algebra we
denote by $(A,Y_A)$, and $\sigma$ extends to an automorphism of $A$ of
order $k$ such that $M$ is a natural weak $\sigma$-twisted $A$-module
in the sense that $Y(\alpha(z),z_1)=\alpha(z_1)$ for $\alpha(z)\in A$
are $\sigma$-twisted vertex operators on $M$.

\begin{rem} {\em $\sigma$ is given by 
$$\sigma a(z)=\lim_{z^{1/k}\to \eta^{-1}z^{1/k}}a(z)$$
 for $a(z)\in A$ (see \cite{Li2}).}
\end{rem}

Let $A^i=\{c(z)\in A| \sigma c(z)=\eta^i c(z)\}$ and $a(z)\in A^{i}$.
For any integer $n$ and $b(z)\in A,$ the operator $a(z)_{n}b(z)$ is an
element of $A$ given by
\begin{eqnarray}\label{3.3}
a(z)_{n}b(z)={\rm Res}_{z_{1}}{\rm
Res}_{z_{0}}\left(\frac{z_{1}-z_{0}}{z}\right)^{i/k}z_{0}^{n}\cdot X
\end{eqnarray}
where 
\[ X=z_{0}^{-1}\delta\left(\frac{z_{1}-z}{z_{0}}\right)a(z_{1})
b(z)-z_{0}^{-1}
\delta\left(\frac{z-z_{1}}{-z_{0}}\right)b(z)a(z_{1}).\]
Or, equivalently, $a(z)_{n}b(z)$ is defined by:
\begin{eqnarray}\label{3.4}
\sum_{n\in \mathbb{Z}}\left(a(z)_{n}b(z)\right)z_{0}^{-n-1}
={\rm Res}_{z_{1}}\left(\frac{z_{1}-z_{0}}{z}\right)^{i/k}
\cdot X.
\end{eqnarray}
Thus following \cite{Li2}, for $a(z)\in A^{i}$, we define $Y_A(a(z),x)$ 
by setting $Y_A(a(z), z_0)b(z)$ equal to (\ref{3.4}).

\begin{lem}\label{l3.5} We have 
\[ [Y_A(Y_g(u^i,z),z_1), Y_A(Y_g(v^j,z),z_2)]=0.\]
\end{lem}

\begin{proof} {}From the vertex algebra structure of $A,$ we have
\begin{eqnarray*}
[Y_A(Y_g(u^i,z),z_1), Y_A(Y_g(v^j,z),z_2)]  \hspace{-1in} \\
& & = \Res_{z_0}z_2^{-1}\delta\left(\frac{z_1-z_0
}{z_2}\right)Y_A(Y_A(Y_g(u^i,z),z_0)Y_g(v^j,z),z_2).
\end{eqnarray*}
So we need to compute $Y_A(Y_g(u^i,z),z_0)Y_g(v^j,z).$

Note that $Y^p_g(u^i,z)=\sum_{n\in p/k+\Z}u^i_nz^{-n-1}$ is an
eigenvector for $\sigma$ with eigenvalue $\eta^p.$ Set
\[ X_p=z_{0}^{-1}\delta\left(\frac{x-z}{z_{0}}\right)Y_g^p(u^i,x)
Y_g(v^j,z)-z_{0}^{-1}\delta\left(\frac{z-x}{-z_{0}}\right)Y_g(v^j,z)
Y_g^p(u^i,x).\]
Then by (\ref{3.4})
\[Y_A(Y_g(u^i,z),z_0)Y_g(v^j,z)=\sum_{p=0}^{k-1}\Res_x\left(\frac{x-z_{0}}{z}\right)^{p/
k}X_p.\] Using Lemma \ref{l3.3} we compute 
\begin{eqnarray*} 
\lefteqn{\Res_{z_0}z_2^{-1}\delta\left(\frac{z_1-z_0
}{z_2}\right)Y_A(Y_g(u^i,z),z_0)Y_g(v^j,z) = } \hspace{.1in}\\
&=& \Res_{z_0}\Res_{x}
\sum_{p=0}^{k-1}z_2^{-1}\delta\left(\frac{z_1-z_0
}{z_2}\right)\left(\frac{x-z_{0}}{z}\right)^{p/k}X_p\\
&=& \Res_{z_0}\Res_{x}
\sum_{p=0}^{k-1}\sum_{n=0}^{\infty}z_2^{-1}\delta\left(\frac{z_1-z_0
}{z_2}\right)\left(\frac{x-z_{0}}{z}\right)^{p/k}\\
& &\hspace{2.1in} z_0^{-1}[Y_g^p(u^i,x),Y_g(v^j,z)](x-z)^nz_0^{-n}\\
&=& \Res_{z_0}\Res_{x}
\sum_{p=0}^{k-1}\sum_{n=0}^{\infty}z_2^{-1}\delta\left(\frac{z_1-z_0
}{z_2}\right)\left(\frac{x-z_{0}}{z}\right)^{p/k} z_0^{-1}(x-z)^nz_0^{-n} \\
& & \hspace{.5in}  \Res_{y}\frac{1}{k}z^{-1}\eta^{(j-i)p}\left(\frac{x-y}
{z}\right)^{-p/k}\delta\left(\frac{x-y}
{z}\right)Y_g((Y(u,y)v)^j,z)\\
&=& \Res_{z_0}\Res_{x}\Res_{y}
\sum_{p=0}^{k-1}\sum_{n=0}^{\infty}z_2^{-1}\delta\left(\frac{z_1-z_0
}{z_2}\right)\left(\frac{x-z_{0}}{z}\right)^{p/k}\\
& & \hspace{.2in}  z_0^{-1}y^nz_0^{-n}
\frac{1}{k}z^{-1}\eta^{(j-i)p}\left(\frac{x-y}
{z}\right)^{-p/k}\delta\left(\frac{x-y}
{z}\right)Y_g((Y(u,y)v)^j,z).
\end{eqnarray*}
Thus using the $\delta$-function identity (\ref{delta function1}), we have
\begin{eqnarray*}
& & \hspace{-1in} \Res_{x}z^{-1}\left(\frac{x-z_{0}}{z}
\right)^{p/k}\left(\frac{x-y}{z}\right)^{-p/k}\delta
\left(\frac{x-y}{z}\right) = \\
&=& \Res_{x}x^{-1}\left(\frac{x-z_{0}}{z}\right)^{p/k}
\left(\frac{z+y}{x}\right)^{p/k}\delta\left(\frac{z+y}{x}\right)\\
&=& \Res_{x}x^{-1}(1-\frac{z_0}{x})^{p/k}z^{-p/k}
(z+y)^{p/k}\delta\left(\frac{z+y}
{x}\right)\\
&=& \Res_{x}x^{-1}(1-\frac{z_0}{z+y})^{p/k}z^{-p/k}
(z+y)^{p/k}\delta\left(\frac{z+y}
{x}\right) \\
&=& z^{-p/k}(z+y-z_0)^{p/k}.
\end{eqnarray*}
Finally, we have
\begin{eqnarray*} 
\lefteqn{\Res_{z_0}z_2^{-1}\delta\left(\frac{z_1-z_0
}{z_2}\right)Y_A(Y_g(u^i,z),z_0)Y_g(v^j,z) =}\\
&=& \Res_{y}\Res_{z_0}\frac{1}{ k}\sum_{p=0}^{k-1}z_2^{-1}\delta
\left(\frac{z_1-z_0}{z_2}\right)\frac{1}{z_0-y}\\
& &\hspace{1.5in} \eta^{(j-i)p}z^{-p/k}(z+y-z_0)^{p/k}
Y_g((Y(u,y)v)^j,z)\\ 
&=& \Res_{y}\Res_{z_0}\frac{1}{ k}\sum_{p=0}^{k-1}z_2^{-1}\delta
\left(\frac{z_1-(z_0+y)}{z_2}\right)\frac{1}{z_0}\\
& &\hspace{1.8in}  \eta^{(j-i)p}z^{-p/k}(z-z_0)^{p/k}
Y_g((Y(u,y)v)^j,z)\\
&=& \Res_{y}\frac{1}{ k}\sum_{p=0}^{k-1}z_2^{-1}\delta\left(\frac{z_1-y}{z_2}
\right)\eta^{(j-i)p}Y_g((Y(u,y)v)^j,z)\\
&=& 0,
\end{eqnarray*}
as desired. \end{proof}

\begin{lem}\label{l3.6} For $u_1,...,u_k\in V$, we have  
\begin{multline*}
Y_A(Y_g(u^k_k,z)_{-1}\cdots Y_g(u_2^2,z)_{-1}Y_g(u_1^1,z),x)  \\
 = Y_A(Y_g(u^k_k,z),x)\cdots Y_A(Y_g(u_2^2,z),x)Y_A(Y_g(u_1^1,z),x)
\end{multline*}
where $Y_g(u_i^i,z)_{-1}$ is the component vertex operator of 
$Y_A(Y(u_i^i,z),x).$
\end{lem}

\begin{proof} The lemma follows {}from Lemma \ref{l3.5} above and 
formula (13.26) of \cite{DL1}.
\end{proof}

Define the map $f : \vtk \rightarrow A$ by
\begin{eqnarray*}
f: \vtk &\rightarrow& A\\
u_1\otimes\cdots \otimes u_k = (u_k^k)_{-1}\cdots (u_2^2)_{-1}u^1_1 &\mapsto&
Y_g(u_k^k,z)_{-1}\cdots 
Y_g(u_2^2,z)_{-1}Y_g(u_1^1,z) 
\end{eqnarray*}
for $u_1,...,u_k\in V$. Then $f(u^i)=Y_g(u^i,z).$ 

\begin{lem}\label{l3.7} $f$ is a homomorphism of vertex algebras.
\end{lem}

\begin{proof} We need to show that
\[fY(u_1\otimes\cdots \otimes u_k,x)=Y_A( Y_g(u_k^k,z)_{-1}\cdots 
Y_g(u_2^2,z)_{-1}Y_g(u_1^1,z),x)f\]
for $u_i\in V.$ 
Take $v_i\in V$ for $i=1,...,k.$ Then 
\begin{eqnarray*}
fY(u_1\otimes\cdots \otimes u_k,x)(v_1\otimes \cdots\otimes v_k) = 
\hspace{-1.4in}\\
&=& \! \! f(Y(u_1,x)v_1\otimes\cdots Y(u_k,x)v_k)\\
&=& \! \! Y_g(Y(u_k^k,x)v_k^k,z)_{-1}\cdots Y_g(Y(u_2^2,x)v_2^2,z)_{-1}
Y_g(Y(u_1^1,x)v_1^1,z) .
\end{eqnarray*}
By Lemma \ref{l3.6}, we have
\[Y_A( Y_g(u_k^k,z)_{-1}\cdots Y_g(u_2^2,z)_{-1}Y_g(u_1^1,z),x)
f(v^1\otimes \cdots \otimes v^k) \hspace{1.5in}\]
\[= Y_A(Y_g(u^k_k,z),x)\cdots Y_A(Y_g(u_2^2,z),x)Y_A(Y_g(u_1^1,z),x)
Y_g(v_k^k,z)_{-1}\]
\[\hspace{3.2in} \cdots Y_g(v_2^2,z)_{-1}Y_g(v_1^1,z).\]
By Lemma \ref{l3.5}, it is enough to show that 
$$ Y_g(Y(u^i,x)v^i,z)=Y_A(Y_g(u^i,z),x)Y_g(v^i,z)$$
for $u,v\in V$ and $i=1,...,k.$ In fact, in view of the
relation between $Y(u^1,z)$ and $Y(u^i,z)$ for 
$u\in V,$  we only need to prove the case $i=1.$ 

By Proposition \ref{psun1},
\begin{eqnarray*}
Y_g(Y(u^1,z_0)v^1,z_2) &=& Y(\Delta_k(z_2)Y(u,z_0)v,z_2^{1/k})\\
&=& Y(Y(\Delta_k(z_2+z_0)u,(z_2+z_0)^{1/k}-z_2^{1/k})\Delta_k(z_2)v,z_2^{1/k}).
\end{eqnarray*}
On the other hand,
\[Y_A(Y_g(u^1,z_2),z_0)Y_g(v^1,z_2)=\sum_{p=0}^{k-1}\Res_{z_1}
\left(\frac{z_1-z_{0}}{z_2}\right)^{p/k}X\]
where
\[X=z_{0}^{-1}\delta\left(\frac{z_1-z_2}{z_{0}}\right)Y_g(u^1,z_1)
Y_g(v^1,z_2)-z_{0}^{-1}\delta\left(\frac{z_2-z_1}{-z_{0}}\right)
Y_g(v^1,z_2)Y_g(u^1,z_1).\]

By equation (\ref{a3.3}), there exists a positive integer $N$ such that
\[(z_1-z_2)^NY_g(u^1,z_1)Y_g(v^1,z_2)=(z_1-z_2)^NY_g(v^1,z_2)Y_g(u^1,z_1).\]
Thus
\begin{eqnarray*}
X &=& z_{0}^{-1}\delta\left(\frac{z_1-z_2}{z_{0}}\right)Y_g(u^1,z_1)
Y_g(v^1,z_2)\\
& &\hspace{1in} - \; z_{0}^{-1} \delta\left(\frac{z_2-z_1}{-z_{0}}\right)
z_0^{-N}(z_1-z_2)^NY_g(v^1,z_2)Y_g(u^1,z_1)\\
&=& z_{0}^{-1}\delta\left(\frac{z_1-z_2}{z_{0}}\right)z_0^{-N}
\left((z_1-z_2)^NY_g(u^1,z_1)Y_g(v^1,z_2)\right)\\
& & \hspace{.8in} - \; z_{0}^{-1} \delta\left(\frac{z_2-z_1}{-z_{0}}\right)
z_0^{-N}\left((z_1-z_2)^NY_g(u^1,z_1)Y_g(v^1,z_2)\right)\\
&=& z_2^{-1}z_0^{-N}\delta\left(\frac{z_1-z_0}{z_2}\right)
\left((z_1-z_2)^NY_g(u^1,z_1)Y_g(v^1,z_2)\right)
\end{eqnarray*}
where we have used the elementary $\delta$-function relation
\begin{equation}\label{a3.4}
z_{0}^{-1}\delta\left(\frac{z_1-z_2}{z_{0}}\right)-z_{0}^{-1}
\delta\left(\frac{z_2-z_1}{-z_{0}}\right)=z_2^{-1}\delta\left(\frac{z_1-z_0}{z_2}\right)
\end{equation}
(cf. \cite{FLM3}). 

Therefore using the $\delta$-function relation (\ref{delta function2}), 
we have
\begin{multline*}
Y_A(Y_g(u^1,z_2),z_0)Y_g(v^1,z_2) \\
= \Res_{z_1} z_0^{-N} z_2^{-1} \delta\Biggl(\frac{(z_1-z_0)^{1/k}}{z_2^{1/k}}\Biggr)
\left((z_1-z_2)^NY_g(u^1,z_1)Y_g(v^1,z_2)\right).
\end{multline*}

Let $x$ be a new formal variable which commutes with $z_0,z_1,z_2.$ 
Then 
\begin{eqnarray*}
& &\hspace{-.5in} z_2^{-1/k}\delta\Biggl(\frac{z_1^{1/k}-x}{z_2^{1/k}}
\Biggr)\left((z_1-z_2)^NY_g(u^1,z_1)Y_g(v^1,z_2)\right) =\\
\hspace{.1in} &=& \! x^{-1}\delta\Biggl(\frac{z_1^{1/k}-z_2^{1/k}}{x}\Biggr)
\left((z_1-z_2)^NY_g(u^1,z_1)Y_g(v^1,z_2)\right)\\
& & \hspace{1in} - \; x^{-1} \delta\Biggl(\frac{-z_2^{1/k}+z_1^{1/k}}{x}
\Biggr)\left((z_1-z_2)^NY_g(u^1,z_1)Y_g(v^1,z_2)\right)\\
&=& \! (z_1-z_2)^Nx^{-1}\delta\Biggl(\frac{z_1^{1/k}-z_2^{1/k}}{x}\Biggr)
Y(\D(z_1)u,z_1^{1/k})Y(\D(z_2)v,z_2^{1/k})\\
& &\hspace{.2in} -\; (z_1-z_2)^Nx^{-1} \delta\Biggl(\frac{-z_2^{1/k} +
z_1^{1/k}}{x}\Biggr)Y(\D(z_2)v,z_2^{1/k})Y(\D(z_1)u,z_1^{1/k})\\
&=& \! (z_1-z_2)^Nz_2^{-1/k}\delta\Biggl(\frac{z_1^{1/k}-x}{z_2^{1/k}}\Biggr)
Y(Y(\D(z_1)u,x)\D(z_2)v,z_2^{1/k}).
\end{eqnarray*}
Note that the first term in the above formula is well defined when
$x$ is replaced by  $z_1^{1/k}-(z_1-z_0)^{1/k}$, and therefore the last 
term is also well defined under this substitution. Thus
\begin{eqnarray*}
& &\hspace{-.35in} z_0^{-N}z_2^{-1/k}\delta\Biggl(
\frac{(z_1-z_0)^{1/k}}{z_2^{1/k}}\Biggr)\left((z_1-z_2)^N
Y_g(u^1,z_1)Y_g(v^1,z_2)\right) =\\
&=& \! z_2^{-1/k}\delta\Biggl(\frac{(z_1-z_0)^{1/k}}{z_2^{1/k}}\Biggr)
Y(Y(\D(z_1)u,z_1^{1/k}-(z_1-z_0)^{1/k})\D(z_2)v,z_2^{1/k})\\
&=& \! z_2^{-1/k}\delta\Biggl(\frac{(z_1-z_0)^{1/k}}{z_2^{1/k}}\Biggr)
Y(Y(\D(z_2+z_0)u,(z_2+z_0)^{1/k}-z_2^{1/k})\D(z_2)v,z_2^{1/k}).
\end{eqnarray*}

Finally we have
\begin{eqnarray*}
& &\hspace{-.4in} Y_A(Y_g(u^1,z_2),z_0)Y_g(v^1,z_2) =\\
&=& \! \Res_{z_1}z_2^{-1}\delta\Biggl(\frac{(z_1-z_0)^{1/k}}{z_2^{1/k}}\Biggr)
Y(Y(\D(z_2+z_0)u,(z_2+z_0)^{1/k}-z_2^{1/k})\\
& & \hspace{3.5in} \D(z_2)v,z_2^{1/k})\\
&=& \! Y(Y(\D(z_2+z_0)u,(z_2+z_0)^{1/k}-z_2^{1/k})\D(z_2)v,z_2^{1/k})\\
&=& \! Y_g(Y(u^1,z_0)v^1,z_2),
\end{eqnarray*}
as desired. \end{proof}

Let $(M,Y)$ be a weak $V$-module.  Define $T_g^k(M,Y) = (T_g^k(M),
Y_g) = (M, Y_g)$.  That is $T_g^k(M, Y)$ is $M$ as the underlying
vector space and the vertex operator $Y_g$ is given by (\ref{define
g-twist}).

Now we state our first main theorem of the paper.
\begin{thm}\label{main1} 
$(T_g^k(M),Y_g)$ is a weak $g$-twisted $V^{\otimes k}$-module such
that $T_g^k(M)=M$, and $Y_g$, defined by (\ref{define g-twist}), is the
linear map {}from $V^{\otimes k}$ to 
\[(\End \; T_g^k(M))[[z^{1/k},z^{-1/k}]]\]
defining the twisted module structure. Moreover, 

(1) $(M, Y)$ is an irreducible weak $V$-module if and only if
$(T_g^k(M), Y_g)$ is an irreducible weak $g$-twisted $V^{\otimes
k}$-module.

(2) $M$ is an admissible $V$-module if and only if $T_g^k(M)$ is an
 admissible $g$-twisted $V^{\otimes k}$-module.

(3) $M$ is an ordinary $V$-module if and only if $T_g^k(M)$ is an
 ordinary $g$-twisted $V^{\otimes k}$-module.
\end{thm}

\begin{proof} It is immediate {}from Lemma \ref{l3.7} that $T_g^k(M)=M$ is a weak
$g$-twisted $V^{\otimes k}$-module with $Y_g(u^1,z)=\bar Y(u,z).$ Note
that
\[Y_g((\Delta_k(z)^{-1}u)^1,z)=\bar Y(\Delta_k(z)^{-1}u,z)=Y(u,z^{1/k})\]
and that all twisted vertex operators $Y_g(v,z)$ for $v\in \vtk$ can
be generated {}from $Y_g(u^1,z)$ for $u\in V.$ It is clear now that $M$
is an irreducible weak $V$-module if and only if $T_g^k(M)$ is an
irreducible weak $g$-twisted $V^{\otimes k}$-module. So (1) has been
proved.

For (2) since $M$ is an admissible $V$-module, we have $M =
\oplus_{n\in \mathbb{N}} M(n)$ such that for $m \in \frac{1}{k}
\mathbb{Z}$, the component operator $u_m$ satisfies $u_mM(n)\subset
M(\wt \; u-m-1+n)$ if $u\in V$ is of homogeneous weight. Define a
$\frac{1}{k}\Z_{+}$-gradation on $T_g^k(M)$ such that $T_g^k(M)(n/k) =
M(n)$ for $n\in \Z.$ Recall that $Y_g(v,z) = \sum_{m \in \frac{1}{k}
\mathbb{Z}} v_nz^{-m-1}$ for $v\in \vtk$.  We have to show that
$v_mT_g^k(M)(n)\subset T_g^k(M)(\wt \; v-m-1+n)$ for $m,n\in\frac{1}{
k}\Z$.  As before, since all twisted vertex operators $Y_g(v,z)$ for
$v\in \vtk$ can be generated {}from $Y_g(u^1,z)$ for $u\in V$, it is
enough to show $u^1_mT_g^k(M)(n)\subset T_g^k(M)(\wt \; u-m-1+n)$.

Let $u\in V_p$ for $p\in\Z$.  Then 
\[\Delta_k(z)u=\sum_{i=0}^{\infty}u(i)z^{1/k - p - i/k}\]
where $u(i)\in V_{p-i}.$ Thus
\[Y_g(u^1,z)=Y(\Delta_k(z)u,z^{1/k})=\sum_{i=0}^{\infty}Y(u(i),z^{1/k}))
z^{1/k - p - i/k} ,\]
and for $m \in \frac{1}{k} \mathbb{Z}$
\[u^1_m=\sum_{i=0}^{\infty}u(i)_{(1-k)p-i-1+km+k}.\]
Since the weight of $u(i)_{(1-k)p-i-1+km+k}$ is $k(p-m-1)$, we see
that $u^1_mT_g^k(M)(n)=u^1_mM(kn)\subset M(k(p-m-1+n))$.  That is
$u^1_mT_g^k(M)(n)\subset T_g^k(M)(p-m-1+n)$, showing that $T_g^k(M)$
is an admissible $g$-twisted $V^{\otimes k}$-module.

Similarly, one can show that if $T_g^k(M)$ is an admissible
$g$-twisted $V^{\otimes k}$-module, then $M$ is an admissible 
$V$-module with $M(n)=T_g^k(M)(n/k)$ for $n\in \Z$. 

In order to prove (3) we write $Y_g(\bar\omega,z) = \sum_{n\in\Z}
L_g(n)z^{-n-2}$ where $\bar\omega=\sum_{j=1}^k\omega^j$.  We have
\[Y_g(\bar\omega,z)=\sum_{i=0}^{k-1}\lim_{z^{1/k}\mapsto \eta^{-i}z^{1/k}}
  Y_g(\omega^1,z).\]
It follows {}from (\ref{sun1}) that 
$L_g(0)=\frac{1}{k}L(0)+\frac{(k^2-1)c}{24k}$.  This immediately implies (3).
\end{proof}

Let $V$ be an arbitrary vertex operator algebra and $g$ an
automorphism of $V$ of finite order. We denote the categories of weak,
admissible and ordinary $g$-twisted $V$-modules by $\mathcal{ C}^g_w(V),$
$\mathcal{ C}^g_a(V)$ and $\mathcal{ C}^g(V)$, respectively.  If $g=1$, we
habitually remove the index $g.$

Now again consider the vertex operator $V^{\otimes k}$ and the
$k$-cycle $g = (1 2 \cdots k)$.  Define
\begin{eqnarray*}
T_g^k: \mathcal{ C}_w(V) &\longrightarrow& \mathcal{ C}^g_w(V^{\otimes k})\\
  (M,Y) &\mapsto& (T_g^k(M),Y_g) = (M,Y_g)\\
     f  &\mapsto& T_g^k(f) = f
\end{eqnarray*}
for $(M,Y)$ an object and $f$ a morphism in $\mathcal{ C}_w(V)$.  

The following corollary to Theorem \ref{main1} is obvious.

\begin{cor}\label{c3.10} 
$T_g^k$ is a functor {}from the category $\mathcal{ C}_w(V)$ to the category
$\mathcal{ C}^g_w(\vtk)$ such that: (1) $T_g^k$ preserves irreducible
objects; (2) The restrictions of $T_g^k$ to $\mathcal{ C}_a(V)$ and $\mathcal{
C}(V)$ are functors {}from $\mathcal{ C}_a(V)$ and $\mathcal{ C}(V)$ to $\mathcal{
C}^g_a(\vtk)$ and $\mathcal{ C}^g(\vtk)$, respectively.
\end{cor}

In the next section we will construct a functor $U_g^k$ {}from the
category $\mathcal{ C}^g_w(\vtk)$ to the category $\mathcal{ C}_w(V)$ such
that $U_g^k \circ T_g^k = id_{\mathcal{ C}_w(V)}$ and $T_g^k \circ U_g^k = 
id_{\mathcal{ C}^g_w(\vtk)}$.

\begin{rem}
In constructing the weak $g$-twisted $V^{\otimes k}$-module structure on 
the weak $V$-module $(M,Y)$, we chose to define $Y_g(u^1,z)$ by 
(\ref{define g-twist}) which then generate $Y_g(v,z)$ for all $v \in 
V^{\otimes k}$.  But our choice of $u^1$ for $u \in V$ as generators is not 
canonical.  We could just as well have chosen to define $Y^j_g(u^j,z) = 
\bar{Y}(u,z)$ for $j = 2,...,k$, and generated $g$-twisted operators 
$Y^j_g(v,z)$ for $v \in \vtk$ {}from these rather than {}from $Y^1_g(u^1,z)$.  
These new $g$-twisted operators $Y^j_g(v,z)$ are related to the old ones $Y_g 
= Y^1_g$ by
\[Y^{j+1}_g(v,z) = \lim_{z^{1/k} \rightarrow \eta^j z^{1/k}} Y_g(v,z) \]
for $j = 1,...,k-1$.  This is a reflection of the fact
that for any vertex operator algebra $(V,Y,\mathbf{1},\omega)$ and any
automorphism $g$ of $V$, if $(M,Y_g)$ is a weak $g$-twisted $V$-module,
then so is $(M, \tilde{Y}_g)$ where \[\tilde{Y}_g(v,z) = \lim_{z^{1/k}
\rightarrow \eta^j z^{1/k}} Y_g(v,z) .\] 
Furthermore, if $(M, Y_g)$ is admissible, then $(M, Y_g)$ and $(M,
\tilde{Y}_g)$ are isomorphic as $g$-twisted $V$-modules via
\begin{eqnarray*}
f : (M, Y_g) & \rightarrow & (M, \tilde{Y}_g) \\
         w & \mapsto & \eta^{jkn}w
\end{eqnarray*}
for $w \in M(n)$.  (Note that here we have used the fact that
$\eta^{k(\mathrm{wt} \; v)} = 1$ for all $v \in V$ of homogeneous weight
$\mathrm{wt} \; v$.)
\end{rem}

\section{Constructing a weak $V$-module structure on a 
weak $g = (12 \cdots k)$-twisted $V^{\otimes k}$-module}
\setcounter{equation}{0}

For $k \in \ZZ$ and $g = (12\cdots k)$, let $M=(M,Y_g)$ be a 
weak $g$-twisted $\vtk$-module.  Motivated by the construction 
of weak $g$-twisted $\vtk$-modules {}from weak $V$-modules in 
Section 3, we consider 
\begin{equation}\label{define U}
Y_g((\Delta_k(z^k)^{-1}u)^1,z^k)
\end{equation}
for $u\in V$ where $\Delta_k (z)^{-1} = \Delta_k^V(z)^{-1}$ is given by 
(\ref{Delta for a module}), i.e.,
\[\Delta_k (z)^{-1} = z^{-\left( 1/k - 1 \right) L(0)}  k^{L(0)} \exp 
\Biggl( -\sum_{j \in \ZZ} a_j z^{- j/k} L(j) \Biggr).  \]
Note that (\ref{define U}) is multivalued since 
$Y_g((\Delta_k(z)^{-1}u)^1,z) \in (\mathrm{End} M) [[z^{1/k}, z^{-1/k}]]$.   
We thus define $Y_U(u,z) = Y_g((\Delta_k(z^k)^{-1}u)^1,z^k)$ to be the  
unique formal Laurent series in $(\mathrm{End} M) [[z, z^{-1}]]$ given by  
taking $(z^k)^{1/k} = z$.  Our goal in this section is to construct a 
functor $U_g^k : \mathcal{ C}_w^g(\vtk) \rightarrow \mathcal{ C}_w(V)$ with
$U_g^k(M,Y) = (U_g^k(M),Y_U) = (M,Y_U)$.  If we instead define $Y_U$ by 
taking $(z^k)^{1/k} = \eta^jz$ for $\eta = e^{-2 \pi i / k}$ with 
$j=1,...,k-1$, then $(M,Y_U)$ will not be a weak $V$-module.  Further 
note that this implies that if we allow $z$ to be complex number and if 
we define $z^{1/k}$ using the principal branch of the logarithm, then much 
of our work in this section is valid if and only if $-\pi/k < 
\mathrm{arg} \; z < \pi/k$.

\begin{lem}\label{l4.1} For $u\in V,$ we have
\begin{eqnarray*}
Y_U(L(-1)u,z) &=& \left(\frac{d}{dz}((z^k)^{1/k})\right)\frac{d}{dz}Y_U(u,z)\\
&=& \frac{d}{dz} Y_U(u,z)
\end{eqnarray*}
on $U_g^k(M) = M$.  Thus the $L(-1)$-derivative property holds for $Y_U$.
\end{lem}

\begin{proof} The proof is similar to that of Lemma \ref{l3.1}.
By Corollary \ref{c2.5} we have
\[\D(z)^{-1}L(-1)-kz^{-1/k + 1}L(-1)\D(z)^{-1} = kz^{-1/k 
+ 1}\frac{d}{dz}\D(z)^{-1}.\]
Making the change of variable $z\to z^k$ gives 
\[\D(z^k)^{-1} L(-1) - k(z^k)^{-1/k} z^k L(-1) \D(z^k)^{-1} =
(z^k)^{-1/k}z\frac{d}{dz}\D(z^k)^{-1}.\]
Thus if $(z^k)^{1/k} = \eta^j z$, we have
\begin{eqnarray*}
& & \hspace{-.6in} \frac{d}{dz}Y_g((\D(z^k)^{-1}u)^1,z^k) =  \\
&=& Y_g((\frac{d}{dz}\D(z^k)^{-1}u)^1,z^k)+\frac{d}{dx}
Y_g((\D(z^k)^{-1}u)^1,x^k)|_{x=z}\\
&=& Y_g((\frac{d}{dz}\D(z^k)^{-1}u)^1,z^k)+
kz^{k-1}Y_g(L(-1)(\D(z^k)^{-1}u)^1,z^k)\\
&=& Y_g((\frac{d}{dz}\D(z^k)^{-1}u)^1,z^k)+kz^{k-1}
Y_g((L(-1)\D(z^k)^{-1}u)^1,z^k)\\
&=& \eta^j Y_g((\D(z^k)^{-1}L(-1)u)^1,z^k).
\end{eqnarray*}
Since by definition $Y_U(u,z) = Y_g((\D(z^k)^{-1}u)^1,z^k)$ with
$(z^k)^{1/k} = z$, the result follows. \end{proof}

\begin{lem}\label{l4.2} Let $u,v\in V$.  Then in $U_g^k(M) = M$, 
\[ [Y_U(u,z_1),Y_U(v,z_2)] \; = \; 
\Res_{z_0} z_2^{-1}\delta\left(\frac{z_1-z_0}{z_2}\right)
Y_U(Y(u, z_0)v,z_2).\] 
\end{lem}

\begin{proof} The proof is similar to the proof of Lemma \ref{l3.2}. {}From the 
twisted Jacobi identity, we have
\begin{equation}\label{commutator for Lemma 4.2}
 [Y_g(u^1,z_1),Y_g(v^1,z_2)] \; = \; \Res_{z_0}\frac{1}{k}z_2^{-1}
\delta\Biggl(\frac{(z_1-z_0)^{1/k}}{z_2^{1/k}}\Biggr)Y_g(
Y(u^1,z_0)v^1,z_2).
\end{equation}
Therefore, 
\begin{eqnarray*}
& & \hspace{-.55in} [Y_U(u,z_1),Y_U(v,z_2)] = \\
&=& [Y_g((\D(z^k_1)^{-1}u)^1,z_1^k),Y_g((\D(z_2^k)^{-1}v)^1,z_2^k)] \\
&=& \Res_{x}\frac{1}{k} z_2^{-k} \delta \left( 
\frac{(z_1^k - x)^{1/k}}{z_2} \right)  Y_g( Y( 
(\D(z_1^k)^{-1}u)^1,x) (\D(z_2^k)^{-1}v)^1,z_2^k).
\end{eqnarray*}
We want to make the change of variable $x=z_1^k-(z_1-z_0)^{k}$ where we
choose $z_0$ such that $((z_1 - z_0)^k)^{1/k} = z_1 - z_0$. Then noting
that $(z_1^k - x)^{n/k} |_{x = z_1^k-(z_1-z_0)^{k}} = (z_1 - z_0)^n$ for 
all $n \in \mathbb{Z}$, and using (\ref{residue change of variables}), 
we have
\begin{eqnarray*}
& &\hspace{-.35in} [Y_U(u,z_1),Y_U(v,z_2)] = \\
&=& \! \! \Res_{z_0} z_2^{-k} (z_1-z_0)^{k-1} \delta \left( \frac{
z_1-z_0}{ z_2} \right) Y_g( Y((\D(z_1^k)^{-1}u)^1,z_1^k -
(z_1-z_0)^{k})\\
& & \hspace{3.45in} (\D(z_2^k)^{-1}v)^1,z_2^k)\\
&=& \! \! \Res_{z_0}  z_2^{-1}\delta\left(\frac{z_1-z_0}{z_2}\right) 
Y_g( Y( (\D(z_1^k)^{-1}u)^1, z_1^k - (z_1-z_0)^{k})(\D(z_2^k)^{-1}v)^1,z_2^k)\\
&=& \! \! \Res_{z_0}  z_2^{-1} \delta \left( \frac{z_1-z_0}{z_2}\right)
Y_g((Y(\D((z_2+z_0)^k)^{-1}u, (z_2+z_0)^k-z_2^{k}) \\
& & \hspace{3.5in} \D(z_2^k)^{-1}v)^1,z_2^k).
\end{eqnarray*}
Thus the proof is reduced to proving
\begin{eqnarray*} 
Y(\D((z_2+z_0)^k)^{-1}u,(z_2+z_0)^k - z_2^{k})\D(z_2^k)^{-1} = \D(z_2^k)^{-1} Y
\left(u, z_0 \right),
\end{eqnarray*} 
i.e., proving
\begin{equation}\label{3.10}
\D(z_2^k)Y(\D((z_2+z_0)^k)^{-1}u,(z_2+z_0)^k-
z_2^{k})\D(z_2^k)^{-1} =  Y \left(u,  z_0 \right).
\end{equation}
In Proposition \ref{psun1}, substituting $u$, $z$ and $z_0$ by
$ \D((z_2+z_0)^k)^{-1}u,$ $z_2^k$ and $(z_2+z_0)^k
- z_2^{k}$, respectively, gives equation (\ref{3.10}).
\end{proof}

\begin{thm}\label{t4.l} With the notations as above, $U_g^k(M,Y_g) = 
(U_g^k(M),Y_U) = (M,Y_U)$ is a weak $V$-module.
\end{thm}

\begin{proof} Since the $L(-1)$-derivation property has been proved for $Y_U$ in 
Lemma \ref{l4.1}, we only need to prove the  Jacobi identity which is
equivalent to the commutator formula given by Lemma  \ref{l4.2} 
and the associator  formula which states that for $u,v\in V$ and $w
\in U_g^k(M)$ there exists  a positive integer $n$ such that
\[(z_0+z_2)^nY_U(u,z_0+z_2)Y_U(v,z_2)w=(z_2+z_0)^nY_U(Y(u,z_0)v,z_2)w .\]

Write $u^1=\sum_{i=0}^{k-1}u^1_{(i)}$ where $gu^1_{(i)} =
\eta^iu^1_{(i)}$.  Then {}from the twisted Jacobi identity, we have the
following associator: there exists a positive integer $m$ such that
for $n\geq m,$
\[(z_0+z_2)^{i/k+n}Y_g(u^1_{(i)},z_0+z_2)Y_g(v^1,z_2)w=
(z_2+z_0)^{i/ k+n}Y_g(Y(u^1_{(i)},z_0)v^1,z_2)w\]
for $i=0,...,k-1$.  Replacing $z_2$ by $z_2^k$ and $z_0$ by 
$(z_0+z_2)^k-z_2^k$ gives 
\begin{multline*}
(z_0+z_2)^{i+kn}Y_g(u^1_{(i)},(z_0+z_2)^k)Y_g(v^1,z_2^k)w\\
= (z_2+z_0)^{i+kn}Y_g(Y(u^1_{(i)},(z_2+z_0)^k-z_2^k)v^1,z_2^k)w.
\end{multline*}
Note that if $a\in \vtk$ such that $ga=\eta^ia$, then $Y_g(a,z)=\sum_{l\in
i/k+\Z}a_nz^{-l-1}$. Thus there exists a positive integer $m_i$
such that if $n_i\geq m_i$, then 
\begin{multline*}
(z_0+z_2)^{n_i}Y_g(u^1_{(i)},(z_0+z_2)^k)Y_g(v^1,z_2^k)w\\
=(z_2+z_0)^{n_i} Y_g(Y(u^1_{(i)},(z_2+z_0)^k-z_2^k)v^1,z_2^k)w
\end{multline*}
for $i=0,...,k-1$.  As a result we see that there exists a positive integer 
$m$ such that if $n\geq m$, then 
\begin{multline*}
(z_0+z_2)^{n}Y_g(u^1,(z_0+z_2)^k)Y_g(v^1,z_2^k)w\\
= (z_2+z_0)^{n}Y_g(Y(u^1,(z_2+z_0)^k-z_2^k)v^1,z_2^k)w.
\end{multline*}

Now write $\Delta_k((z_0+z_2)^k)^{-1}u= \sum_{j \in \mathbb{N}} 
u_j(z_0+z_2)^{s_j}$ for some $u_j\in V$ and integers $s_j$, and note that 
this is a finite sum.  Similarly we have a finite sum
$\Delta_k(z_2^k)^{-1}v=\sum_{j\in \mathbb{N}}v_jz_2^{t_j}$ for some $v_j \in
V$ and $t_j \in \mathbb{Z}$.  Thus there exists a positive integer
$m$ such that if $n\geq m$, then
\begin{multline*}
(z_0+z_2)^{n+s_i}Y_g(u_i^1,(z_0+z_2)^k)Y_g(v_j^1,z_2^k)w\\
=(z_2+z_0)^{n+s_i}Y_g(Y(u_i^1,(z_2+z_0)^k-z_2^k)v_j^1,z_2^k)w
\end{multline*}
for all $i,j\in \mathbb{N}$.  Finally we have for $n\geq m,$ 
\begin{eqnarray*}
& &\hspace{-.45in}(z_0+z_2)^nY_U(u,z_0+z_2)Y_U(v,z_2)w = \\
&=& \! (z_2+z_0)^nY_g((\D((z_0+z_2)^k)^{-1}u)^1,
(z_0+z_2)^k)Y_g((\D(z_2^k)^{-1}v)^1,z_2^k) w\\
&=& \! \sum_{i,j \geq 0}(z_0+z_2)^{n+s_i}z_2^{t_j}Y_g(u_i^1,(z_0+z_2)^k)
Y_g(v_j^1,z_2^k)w\\
&=& \! \sum_{i,j \geq 0}(z_2+z_0)^{n+s_i}z_2^{t_j}Y_g(Y(u_i^1,
(z_2+z_0)^k-z_2^k)v_j^1,z_2^k)w\\
&=& \! (z_2+z_0)^nY_g(Y(\D((z_2+z_0)^k)^{-1}u)^{1},(z_2+z_0)^k-z_2^k)
(\D(z_2^k)^{-1}v)^1,z_2^k)w\\
&=& \! (z_2+z_0)^nY_g((\D(z_2^k)^{-1}Y(u,z_0)v)^1,z_2^k)w\\
&=& \! (z_2+z_0)^nY_U(Y(u,z_0)v,z_2)w
\end{eqnarray*}
where we have used equation (\ref{3.10}), completing the proof.
\end{proof}

\begin{thm}\label{t4.ll} 
$U_g^k$ is a functor {}from the category $\mathcal{ C}_w^g(\vtk)$ of
weak $g$-twisted $\vtk$-modules to the category $\mathcal{ C}_w(V)$ of
weak $V$-modules such that $T_g^k \circ U_g^k = 
id_{\mathcal{ C}_w^g(\vtk)}$ and $U_g^k \circ T_g^k = 
id_{\mathcal{ C}_w(V)}$.  In particular, the categories 
$\mathcal{ C}_w^g(\vtk)$ and $\mathcal{ C}_w(V)$ are isomorphic. Moreover,

(1) The restrictions of $T_g^k$ and $U_g^k$ to the category of
admissible $V$-modules $\mathcal{ C}_a(V)$ and to the category of
admissible $g$-twisted $\vtk$-modules $\mathcal{ C}_a^g(\vtk)$,
respectively, give category isomorphisms. In particular, $V$ is
rational if and only if $\vtk$ is $g$-rational.

(2) The restrictions of $T_g^k$ and $U_g^k$ to the category of
ordinary $V$-modules $\mathcal{ C}(V)$ and to the category of ordinary
$g$-twisted $\vtk$-modules $\mathcal{ C}^g(\vtk)$, respectively, give
category isomorphisms.
\end{thm}

\begin{proof} It is trivial to verify $T_g^k \circ U_g^k = id_{\mathcal{
C}_w^g(\vtk)}$ and $U_g^k \circ T_g^k = id_{\mathcal{ C}_w(V)}$ {}from the
definitions of the functors $T_g^k$ and $U_g^k$.  Parts 1 and 2 follow {}from Theorem \ref{main1}.  \end{proof}

\section{Twisted sectors for an arbitrary $k$-cycle}
\setcounter{equation}{0}

In this section we extend Theorem 4.4 to the category of weak
$g$-twisted $V^{\otimes k}$-modules where $g$ is an arbitrary
$k$-cycle.  But we first recall some general results following, for
example, \cite{DLM3}.  

Let $V$ be any vertex operator algebra, and let $g, h \in
\mathrm{Aut}(V)$.  There is an isomorphism {}from the category of weak
$g$-twisted $V$-modules $\mathcal{ C}_w^g(V)$ to the category of weak
$hgh^{-1}$-twisted $V$-modules $\mathcal{ C}_w^{hgh^{-1}}(V)$ given by
\begin{eqnarray*}
h: \mathcal{ C}_w^g(V) &\longrightarrow& \mathcal{ C}^{hgh^{-1}}_w(V)\\
  (M,Y_g) &\mapsto& (M, Y_{hgh^{-1}})
\end{eqnarray*}
where
\[Y_{hgh^{-1}} (u,z) = Y_g (hu,z) \]
for $u \in V$.  Moreover it is clear that $(M, Y_g)$ is irreducible, 
admissible, or ordinary if and only if $(M, Y_{hgh^{-1}})$ is 
irreducible, admissible, or ordinary, respectively.

Now suppose $g'$ is an arbitrary $k$-cycle and let $g = (12 \cdots
k)$.  Then there exists $h \in S_k$ such that $g' = hgh^{-1}$.
However, this $h$ is unique only up to multiplication on the right by
powers of $g$.  Thus given $g'$, we can specify $h$ uniquely by
requiring that $h$ leave 1 fixed.  Denote this by $h_1$.  Then we have
the following Corollary to Theorem 4.4.

\begin{cor}\label{arbitrary k-cycle}
Let $g' \in S_k$ be a $k$-cycle, let $g = (12 \cdots k)$, and let $h_1$
be the unique element of $S_k$ that fixes 1 and satisfies $g' =
h_1gh_1^{-1}$.  Then we have the following isomorphism of categories
\begin{eqnarray*}
T_{g'}^k = h_1 \circ T_g^k : \mathcal{ C}_w(V) &\longrightarrow& \mathcal{
C}^{g'}_w(V^{\otimes k})\\
(M,Y) &\mapsto& (T_{g'}^k(M), Y_{g'}) = (M, Y_{h_1gh_1^{-1}})
\end{eqnarray*}
where 
\[Y_{g'}(v,z) = Y_{h_1gh_1^{-1}} (v,z) = Y_g (h_1v,z) \]
for $v \in V^{\otimes k}$, and $Y_g$ is uniquely determined by
\[Y_g (u^1,z) = \bar{Y}(u,z) = Y(\Delta_k(z)u, z^{1/k})\]
\[Y_g(u^{j+1},z) = \lim_{z^{1/k} \rightarrow \eta^{-j}z^{1/k}} 
Y_g(u^1,z)\]
for $u \in V$.  Moreover, $T_{g'}^k$ preserves irreducible, admissible 
and ordinary objects. 
\end{cor}

\begin{rem} 
In the corollary above, we could just have easily 
defined $T_{g'}^k$ by $T_{g'}^k = h_l \circ T_g^k$ where $h_l$, rather
than fixing 1, is the unique element of $S_k$ that takes 1 to $l$ for
$l = 2,...,k$ and satisfies $g' = h_lgh_l^{-1}$.  Then these new
$g'$-twisted operators defined on $M$ denoted by $Y_{g'}^l$ would
differ {}from those defined in Corollary \ref{arbitrary k-cycle},
denoted by $Y_{g'}$, by 
\[Y^{l + 1}_{g'}(u^1,z) = \lim_{z^{1/k} \rightarrow \eta^l z^{1/k}} 
Y_{g'}(u^1,z)\]
for $l = 1,...,k - 1$ and $u \in V$, where of course $Y_{g'}^{l+1}(v,z)$
for $v \in \vtk$ is determined by $Y_{g'}^{l+1}(u^{j+1},z) = \lim_{z^{1/k} 
\rightarrow \eta^{-j}z^{1/k}}  Y_{g'}^{l+1}(u^1,z)$. 
\end{rem}

\section{Twisted sectors for arbitrary permutations}
\setcounter{equation}{0}

In order to determine the various twisted module categories for an
arbitrary permutation, we need to study twisted modules for tensor
product vertex operator algebras in general.

For $s \in \ZZ$ and $i=1,...s$, let $V^i=(V^i,Y_i,{\bf 1}_i,\omega_i)$ be a 
vertex operator algebra and $g_i$ an automorphism of $V^i$ of finite
order. Then $g=g_1\otimes \cdots \otimes g_s$ is an automorphism of
the tensor product vertex operator algebra $V^1\otimes \cdots\otimes
V^s$ of finite order.  If $W^i$ are weak $g_i$-twisted $V^i$-modules
for $i=1,...,s$ then the tensor product $W^1\otimes \cdots\otimes W^s
= W$ is a weak $g$-twisted $V^1\otimes \cdots\otimes V^s$-module in
the obvious way. Thus {}from Proposition 4.7.2 in \cite{FHL}, we have the
following lemma.

\begin{lem}\label{l5.1} If $W^i$ are irreducible weak 
$g_i$-twisted $V^i$-modules for $i=1,...,s$, then $W = W^1 \otimes 
\cdots \otimes W^s$ is an irreducible
weak $g$-twisted $V^1\otimes \cdots\otimes V^s$-module.
\end{lem}

\begin{prop}\label{p5.2} The notation is the same as before. If 
$V^i$ is $g_i$-rational for all $i$ then $V^1\otimes \cdots \otimes 
V^s$ is $g$-rational and  each irreducible $V^1\otimes \cdots \otimes 
V^s$-module $W$ is isomorphic to $W^1\otimes \cdots\otimes W^s$ for 
some irreducible $g_i$-twisted $V^i$-modules $W^i$.
\end{prop}

\begin{proof} Let $W=\oplus_{n\in\Q_{+}}W(n)$ be an admissible $g$-twisted
$V^1\otimes \cdots\otimes V^s$-module.  We need to show that $W$ is
completely reducible. It is sufficient to prove that a submodule
generated by any vector $w\in W(n)$ is completely reducible, and thus
we can assume that $W$ is generated by $w$.  Identify $V^i$ with the 
vertex operator algebra ${\bf 1}_1\otimes \cdots {\bf 1}_{i-1}\otimes
V^i\otimes {\bf 1}_{i+1}\otimes \cdots \otimes {\bf 1}_s$.  (This is
``almost'' a subalgebra of $V^1\otimes \cdots\otimes V^s$, failing to
be a subalgebra due to the fact that ${\bf 1}_1 \otimes \cdots 
{\bf 1}_{i-1} \otimes V^i\otimes {\bf 1}_{i+1}\otimes \cdots \otimes 
{\bf 1}_s$ has a different Virasoro element {}from that of $V^1 \otimes
\cdots \otimes V^s$.)  With this identification, the submodule $W$ 
generated by $w$ is an admissible $g_i$-twisted $V^i$-module.  By 
\cite{DM} and \cite{Li1}, $W$ is spanned by $\{v^1_{m_1}\cdots v^s_{m_s}w|v^i\in 
V^i, m_i\in \Q\}$.  Set $W^i$ to be the span of $\{v_{m}w|v\in V^i, m
\in \Q\}$.  Since $V^i$ is $g_i$-rational, $W^i$ is an admissible 
$g_i$-twisted $V^i$-module which is a direct sum of irreducible 
$g_i$-twisted $V^i$-modules.  Note that the map $f$ {}from $W^1\otimes 
\cdots\otimes W^s$ to $W$ which sends $v^1_{m_1}w\otimes \cdots 
\otimes v^s_{m_s}w$ to $v^1_{m_1}\cdots v^s_{m_s}w$ for $v^i\in V^i, 
m_i\in \Q$ is a $V^1\otimes \cdots\otimes V^s$-homomorphism.  By Lemma 
\ref{l5.1}, $W^1\otimes \cdots\otimes W^s$ is completely reducible.  
Thus $W$ as a homomorphic image of $W^1\otimes \cdots\otimes W^s$ is 
completely reducible.  In particular, if $W$ is irreducible then $W$ 
is isomorphic to a tensor product of irreducible $g_i$-twisted
$V^i$-modules. \end{proof}

\begin{rem} 
Proposition \ref{p5.2} is a generalization of a special case of a
result proved in \cite{FHL} which states that in the case $g_i=1$ for 
all $i$, any ordinary irreducible $V^1\otimes \cdots\otimes V^s$-module 
$W$ on whose lowest weight space the operators $L_i(0)$ have only 
rational eigenvalues is isomorphic to the tensor product of some 
irreducible $V^i$-modules.  In our case $L_i(0)$ is the component 
operator of $Y_i(\o_i,z)$ of weight 0, the $g_i$ are general 
automorphisms of $V^i$ of finite order, but we have restricted to
the case where each $V^i$ is $g_i$-rational. 
\end{rem}

Let $k$ be a fixed positive integer, and let $g\in S_k$.  Then $g$ can be
written as a product of disjoint cycles $g=g_1\cdots g_p$ where the
order of $g_i$ is $k_i$ such that $\sum_{i}k_i=k$.  (Note that we are
including 1-cycles.)  Furthermore, there exists $h \in S_k$ satisfying
$g=hg_1'\cdots g_p'h^{-1}$ such that $g_i'$ is a $k_i$-cycle which
permutes the numbers $(\sum_{j = 1}^{i-1} k_j) + 1,(\sum_{j = 1}^{i-1} 
k_j) + 2, ..., \sum_{j = 1}^i k_j$.  In Sections 3, 4 and 5 we have 
determined various $g_i$-twisted $\vti$-module categories in terms of 
corresponding $V$-module categories via the functor $T_{g_i}^{k_i}$ 
for $g_i$ a $k_i$-cycle. Thus we have the following:

\begin{thm}\label{t5.3}  Let $W_i$ be a weak $V$-module for $i=1,...,p$.  
Given $g \in S_k$ and a decomposition $g = hg_1'\cdots g_p'h^{-1}$ as
above, $h \circ (T^{k_1}_{g_1'}(W_1)\otimes \cdots \otimes
T^{k_p}_{g_p'}(W_p))$ is a weak $g$-twisted $\vtk$-module.  Moreover, 

(1) If each $W_i$ is irreducible then $h \circ (T^{k_1}_{g_1'}(W_1)
\otimes \cdots\otimes T^{k_p}_{g_p'}(W_p))$ is irreducible.  

(2) If each $W_i$ is an admissible (resp., ordinary) $V$-module, then 
$h \circ (T^{k_1}_{g_1'}(W_1)\otimes \cdots \otimes T^{k_p}_{g_p'}
(W_p))$ is an admissible (resp., ordinary) $g$-twisted $\vtk$-module.

(3) If $V$ is a rational vertex operator algebra, then $\vtk$ is
$g$-rational and all irreducible $g$-twisted modules are given by $h
\circ (T^{k_1}_{g_1'}(W_1)\otimes \cdots \otimes T^{k_p}_{g_p'}(W_p))$
where $W_i$ are irreducible $V$-modules.
\end{thm}

\begin{proof} Parts (1) and (2) follow {}from Corollary \ref{arbitrary k-cycle} 
and Lemma \ref{l5.1}.  Part (3) follows {}from Corollary 
\ref{arbitrary k-cycle} and Proposition \ref{p5.2}. \end{proof}

\begin{rem} 
Note that in the Theorem above, we have not specified 
a unique decomposition $g = hg_1'\cdots g_p'h^{-1}$ but rather have
given the functor $h \circ (T^{k_1}_{g_1'} \otimes \cdots \otimes
T^{k_p}_{g_p'})$ for a given such (non-unique) decomposition 
$g = hg_1'\cdots g_p'h^{-1}$.  However, for any decomposition $g = 
hg_1'\cdots g_p'h^{-1}$, the resulting $g$-twisted $\vtk$-module $h 
\circ (T^{k_1}_{g_1'}(W_1)\otimes \cdots \otimes T^{k_p}_{g_p'}(W_p))$ 
is isomorphic to that obtained {}from any other decomposition.  Similar 
to the symmetries discussed in Remarks 3.11 and 5.2 above, the
isomorphism consists of transformations $\lim_{z^{1/k_i} \rightarrow
\eta^{-j_i}_{k_i} z^{1/k_i}}$ of each $i$-th tensor factor in the
$g$-twisted vertex operator where $\eta_{k_i} = e^{-2\pi i/k_i}$ and 
$j_i = 1,..., k_i - 1$. 
\end{rem}

\section{Twisted sectors for the product of a permutation
with an automorphism of $V$}

In this section we give a slight generalization of our previous results
to a broader class of automorphisms. Note that $\Aut (V)$ acts on 
$\vtk$ diagonally and this action commutes with the action of $S_k$.  
Let $\gamma \in \Aut (V)$ and $g\in S_k$.  In this section we determine 
various $\gamma g$-twisted $\vtk$-module categories.  (Although the 
diagonal action of $\gamma$ on $\vtk$ is more appropriately denoted by
$\gamma^{\otimes k}$, to simplify notation, we will write this as 
$\gamma$ with the diagonal action being understood.)

First, we assume that $g=(12\cdots k)$ and show that suitable variants 
of all the arguments and results in Sections 3 and 4 remain valid in 
for the case of $\gamma g$-twisted $\vtk$-module categories.  Then we 
generalize the permutation $g$ as we did in Sections 5 and 6 and show
that again suitable variants of our arguments remain valid. 

We proceed to formulate these results.  Let $g = (12\cdots k)$. For each
positive integer $n$, set $\eta_n=e^{-2\pi i/n}.$ Then
$\eta_{n/m}=\eta_n^m.$ Let $o(\gamma)=l$, and set $o(\gamma) = l =sd$ 
where $d$ is the greatest common divisor of $l$ and $k$.  Furthermore, 
if we let $o(\gamma g) = n$, then $o(\gamma g)=n=sk$.  Let $(W,Y)$ be a 
weak $\gamma g$-twisted $\vtk$-module.  Then the analogue of equation 
(\ref{k1}) is
\[z^{-1}_0\delta\left(\frac{z_1-z_2}{z_0}\right)Y(u^1,z_1)Y(v^1,z_2)-
z^{-1}_0\delta\left(\frac{z_2-z_1}{-z_0}\right)Y(v^1,z_2)Y(u^1,z_1)\]
\begin{equation}\label{k3}
=\frac{1}{n}z_2^{-1}\sum_{j=0}^{n-1}\delta\Biggl(\eta^j_n\frac{(z_1-z_0)^
{1/n}}{z_2^{1/n}}\Biggr)Y(Y(\gamma^jg^ju^1,z_0)v^1,z_2)
\end{equation} 
for $u,v\in V.$ Let $u$ be an eigenvector for $\gamma$ with eigenvalue
$\eta_{l}^r$ for some $r = 1,...,l$.  Then $\gamma^jg^ju^1=\eta_{l}^{rj}
u^{j+1}$ where $u^{j+k}=u^j$.  Thus similar to the situation in 
Section 3, $Y(\gamma^jg^ju^1,z_0)v^1$ only involves nonnegative integer 
powers of $z_0$ if $j\ne 0\ (\mod\ k)$, and
\begin{eqnarray}\label{k4}
& &\hspace{-.65in} [Y(u^1,z_1),Y(v^1,z_2)] = \nonumber\\
&=& \! \sum_{q=0}^{s-1}\Res_{z_0}\frac{1}{n}z_2^{-1}\delta\Biggl(\eta_n^{kq}
\frac{(z_1-z_0)^{1/n}}{z_2^{1/n}}\Biggr)Y(
Y(\eta_{l}^{rkq}u^1,z_0)v^1,z_2)\nonumber\\
&=& \! \sum_{q=0}^{s-1}\Res_{z_0}\frac{1}{n}z_2^{-1}\delta \Biggl(\eta_s^{q}
\frac{(z_1-z_0)^{1/n}}{z_2^{1/n}}\Biggr)\eta_{s}^{rqk/d}Y(
Y(u^1,z_0)v^1,z_2)\nonumber\\
&=& \! \sum_{q=0}^{s-1}\sum_{p=0}^{n-1}\Res_{z_0}\frac{1}{n}z_2^{-1}
\eta_{s}^{qp}  \left(\frac{z_1-z_0}{z_2}\right)^{p/n}  \delta\left(
\frac{z_1-z_0}{z_2}\right)  \eta_{s}^{rqk/d} \nonumber\\
& & \hspace{2.8in} Y(Y(u^1,z_0)v^1,z_2)\nonumber\\
&=& \! \sum_{q=0}^{s-1}\sum_{t=1}^{k}\Res_{z_0}\frac{1}{sk}z_2^{-1} 
\left(\frac{z_1-z_0}{z_2}\right)^{(ts-rk/d)/n}  \delta\left(
\frac{z_1-z_0}{z_2}\right)  \nonumber\\
& & \hspace{2.8in} Y(Y(u^1,z_0)v^1,z_2)\nonumber\\
&=& \! \sum_{t=1}^{k}\Res_{z_0}\frac{1}{k}z_2^{-1} \left(\frac{z_1-z_0}{z_2}
\right)^{t/k}  \left(\frac{z_1-z_0}{z_2}\right)^{-r/l}  \delta\left(
\frac{z_1-z_0}{z_2}\right) \nonumber\\
& & \hspace{2.8in}  Y(Y(u^1,z_0)v^1,z_2)\nonumber\\
&=& \! \Res_{z_0}\frac{1}{k}z_2^{-1}\left(\frac{z_1-z_0}{z_2}\right)^{-r/l}
\delta\Biggl(\frac{(z_1-z_0)^{1/k}}{z_2^{1/k}}\Biggr)Y(Y(u^1,z_0)v^1,z_2).
\end{eqnarray}
As in Section 3 this shows that for $u\in V$ the component operators of 
$Y(u^1,z)$ on $W$ form a Lie algebra.

Let $(M,Y)$ be a weak $\gamma^k$-twisted $V$-module. Set $\bar
Y(u,z)=Y(\D(z)u,z^{1/k})$.

\begin{lem}\label{l6.1} The assertion of Lemma \ref{l3.1} holds in the 
present setting: for $u\in V$
\[\bY(L(-1)u,z)=\frac{d}{dz}\bY(u,z).\]
\end{lem}

\begin{proof} The steps of the proof for Lemma \ref{l3.1} still hold in the present
setting for $(M,Y)$ a $\gamma^k$-twisted $V$-module rather than just a 
$V$-module. \end{proof}

The analogue of Lemma \ref{l3.2} is  

\begin{lem}\label{l6.2} For $u,v\in V$ such that $\gamma u=\eta_{l}^ru$, 
we have 
\begin{equation}\label{commutator for gamma}
[\bY(u,z_1),\bY(v,z_2)]  =  \Res_{z_0}\frac{1}{k}z_2^{-1} 
\left(\frac{z_1-z_0}{z_2}\right)^{-r/l} \! \delta \Biggl(
\frac{(z_1-z_0)^{1/k}}{z_2^{1/k}}\Biggr)\bY(Y(u,z_0)v,z_2).
\end{equation}
\end{lem}

\begin{proof} {}From the twisted Jacobi identity, we have
\begin{equation}\label{commutator for Lemma 7.2}
 [Y(u,z_1),Y(v,z_2)]  =  \Res_x z_2^{-1} 
\left(\frac{z_1-x}{z_2}\right)^{-r/l} \delta \Biggl(
\frac{z_1-x}{z_2}\Biggr)Y(Y(u,x)v,z_2) 
\end{equation}
for $u \in V$ such that $\gamma u = \eta_{l}^ru$.  Following the steps of 
the proof of Lemma \ref{l3.2} using the commutator formula 
(\ref{commutator for Lemma 7.2}) instead of 
(\ref{commutator for Lemma 3.3}) gives (\ref{commutator for gamma}).
\end{proof}

For $u\in V$ such that $\gamma u=\eta_{l}^r$ define operators $Y_{\gamma 
g}$ on $\vtk$ by
\begin{eqnarray}
Y_{\gamma g}(u^1,z) = \bY(u,z) = Y(\Delta_k(z)u,z^{1/k}) \label{Y for gamma g} \\
\eta_{l}^{rj}Y_{\gamma g}(u^{j+1},z)=\lim_{z^{1/n}\to \eta_n^{-j}z^{1/n}}
Y(u^1,z) \label{Y for gamma g 2}
\end{eqnarray}
for $j = 0,...,k-1$.  Then as before, the operators $Y_{\gamma g}(u^j,z)$
for $u\in V$ and $j=1,...,k$ are mutually local and generate a local
system $A$ in the sense of \cite{Li2}. Let $\sigma$ be a map {}from $A$ to $A$
such that $\sigma Y_{\gamma g}(u^j,z)=\eta_{l}^{-r}Y_{\gamma g}(u^{j+1},
z)$ for $u\in V$ satisfying $\gamma u = \eta^r_l u$, and for $j=1,...,k$.  
Again $A$ has the structure of a vertex algebra $(A,Y_A)$ and $\sigma$ 
extends to an automorphism of $A$ of order $n$ such that $M$ is naturally
a weak $\sigma$-twisted $A$-module in the sense that
$Y(\alpha(z),z_1)=\alpha(z_1)$ for $\alpha(z)\in A.$

\begin{lem}\label{l6.3} The assertions of Lemmas \ref{l3.5}, \ref{l3.6}
and \ref{l3.7} hold for $Y_{\gamma g}$ and the corresponding local
system.
\end{lem}

\begin{proof} The steps of the proofs of Lemmas \ref{l3.5}, \ref{l3.6}
and \ref{l3.7} remain valid in the present setting. \end{proof}

Let $(M,Y)$ be a weak $\gamma^k$-twisted $V$-module.  Define 
$$T_{\gamma g}^k (M,Y)=(T_{\gamma g}^k(M),Y_{\gamma g}) = (M,Y_{\gamma g}).$$
That 
is $T_{\gamma g}^k(M,Y)$ is $M$ as the underlying vector space and the 
vertex operator $Y_{\gamma g}$ is given by (\ref{Y for gamma g}) and 
(\ref{Y for gamma g 2}) .

Now repeating the proof of Theorem \ref{main1} under the present
circumstances, we have the following result which is a generalization 
of Theorem \ref{main1}.
 
\begin{thm}\label{t6.4} 
$(T_{\gamma g}^k(M),Y_{\gamma g})$ is a weak $\gamma g$-twisted 
$V^{\otimes k}$-module such that $T_{\gamma g}^k(M)\!=M$ and 
$Y_{\gamma g}$, defined by (\ref{Y for gamma g}) and (\ref{Y for gamma 
g 2}), is the linear map {}from $V^{\otimes n}$ to 
\[(\End \; T_{\gamma g}^k(M))[[z^{1/n},z^{-1/n}]]\] 
defining the twisted module structure where $n = o(\gamma g)$.  Moreover,

(1) $(M,Y)$ is an irreducible weak $\gamma^k$-twisted $V$-module if and 
only if $(T_{\gamma g}^k(M), Y_{\gamma g})$ is an irreducible weak 
$\gamma g$-twisted $V^{\otimes k}$-module.

(2) $M$ is an admissible $\gamma^k$-twisted $V$-module if and only if
 $T_{\gamma g}^k(M)$ is an admissible $\gamma g$-twisted $V^{\otimes 
k}$-module.

(3) $M$ is an ordinary $\gamma^k$-twisted $V$-module if and only if
 $T_{\gamma g}^k(M)$ is an ordinary $\gamma g$-twisted $V^{\otimes 
k}$-module.
\end{thm}

Let
\begin{eqnarray*}
T_{\gamma g}^k: \mathcal{ C}^{\gamma^k}_w(V) &\longrightarrow& \mathcal{
  C}^{\gamma g}_w(V^{\otimes k})\\ 
  (M,Y) &\mapsto& (T_{\gamma g}^k(M),Y_{\gamma g}) = (M,Y_{\gamma g})\\
  f &\mapsto& T_{\gamma g}^k(f) = f
\end{eqnarray*}
for $(M,Y)$ an object and $f$ a morphism in $\mathcal{ C}^{\gamma^k}_w(V)$.
We have the following corollary to Theorem \ref{t6.4}.

\begin{cor}\label{c6.5} $T_{\gamma g}^k$ is a functor {}from the category 
$\mathcal{ C}_w^{\gamma^k}(V)$ to the category $\mathcal{ C}^{\gamma g}_w
(\vtk)$ such that: (1) $T_{\gamma g}^k$ preserves irreducible objects; 
(2) The restrictions of $T_{\gamma g}^k$ to $\mathcal{ C}_a^{\gamma^k}(V)$ 
and $\mathcal{ C}^{\gamma^k}(V)$ are functors {}from $\mathcal{ C}_a^{\gamma^k}
(V)$ and $\mathcal{ C}^{\gamma^k}(V)$ to $\mathcal{ C}^{\gamma g}_a(\vtk)$ and 
$\mathcal{ C}^{\gamma g}(\vtk)$, respectively.
\end{cor}

Note that if $\gamma$ is the identity, Theorem \ref{t6.4} and Corollary
\ref{c6.5} reduce to Theorem \ref{main1} and Corollary \ref{c3.10},
respectively.

Now let $M=(M,Y_{\gamma g})$ be a weak $\gamma g$-twisted $\vtk$-module.  
In analogy to Section 4, we set $U_{\gamma g}^k(M)=M$ and define
$Y_U(u,z)=Y_{\gamma g}((\Delta_k(z^k)^{-1}u)^1,z^k)$ for $u\in V$ to be 
the unique power series obtained by letting $(z^k)^{1/k} = z$.  Lemma
\ref{l4.1} remains valid in this case and Lemma \ref{l4.2} is modified
as follows:

\begin{lem}\label{l6.7} 
 Let $u,v\in V$ such that $\gamma u=\eta_{l}^r.$  Then on $U_{\gamma 
g}^k(M)$, 
\begin{equation}\label{U commutator for gamma} 
[Y_U(u,z_1),Y_U(v,z_2)]  =  
\Res_{z_0}z_2^{-1} \left( \frac{z_1-z_0}{z_2} \right)^{-rk/l} \!
\delta\left(\frac{z_1-z_0}{z_2}\right) Y_U(Y(u,z_0)v,z_2). 
\end{equation}
\end{lem}

\begin{proof}  The $\gamma g$-twisted vertex operators satisfy
\begin{multline}\label{commutator for Lemma 7.6}
[Y_{\gamma g}(u^1,z_1),  Y_{\gamma g}(v^1,z_2)]   \\ 
=\Res_{z_0} \frac{1}{k} z_2^{-1} \left( \frac{z_1-z_0}{z_2} \right)^{-r/l} 
\delta\Biggl(\frac{(z_1-z_0)^{1/k}}{z_2^{1/k}}\Biggr) Y_{\gamma g}
(Y(u^1,z_0)v^1,z_2). 
\end{multline}
Following the steps for the proof of Lemma \ref{l4.2} using
(\ref{commutator for Lemma 7.6}) instead of (\ref{commutator for Lemma 4.2}) 
gives (\ref{U commutator for gamma}). 
\end{proof}

Repeating the proof of Theorem \ref{t4.l} in the present setting gives:
 
\begin{thm}\label{t6.8} With the notations as above, $(U_{\gamma 
g}^k(M),Y_U)$ is a weak $\gamma^k$-twisted $V$-module.
\end{thm} 

Thus we have the following analogue of Theorem \ref{t4.ll}.

\begin{thm} 
$U_{\gamma g}^k$ is a functor {}from the category $\mathcal{ C}_w^{\gamma 
g}(\vtk)$ of weak $\gamma g$-twisted $\vtk$-modules to the category 
$\mathcal{ C}^{\gamma^k}_w(V)$ of weak $\gamma^k$-twisted $V$-modules 
such that $T_{\gamma g}^k \circ U_{\gamma g}^k = 
id_{\mathcal{ C}_w^{\gamma g}(\vtk)}$ and $U_{\gamma g}^k \circ  
T_{\gamma g}^k = id_{\mathcal{ C}_w^{\gamma^k}(V)}$.  In particular, 
the categories $\mathcal{ C}_w^{\gamma g}(\vtk)$ and 
$\mathcal{ C}_w^{\gamma^k}(V)$ are isomorphic.  Moreover,

(1) The restrictions of $T_{\gamma g}^k$ and $U_{\gamma g}^k$ to the 
category of admissible $V$-modules $\mathcal{ C}_a^{\gamma^k}(V)$ and to 
the category of admissible $\gamma g$-twisted $\vtk$-modules 
$\mathcal{ C}_a^{\gamma g}(\vtk)$, respectively, give category 
isomorphisms. In particular, $V$ is $\gamma^k$-rational if and only 
if $\vtk$ is $\gamma g$-rational.

(2) The restrictions of $T_{\gamma g}^k$ and $U_{\gamma g}^k$ to the 
category of ordinary $V$-modules $\mathcal{ C}^{\gamma^k}(V)$ and to the 
category of ordinary $\gamma g$-twisted $\vtk$-modules 
$\mathcal{ C}^{\gamma g}(\vtk)$, respectively, give category isomorphisms.
\end{thm}

Recalling the conjugation functor $h$ {}from Section 5, and again noting
that the diagonal action of $\gamma$ on $V^{\otimes k}$ commutes with 
that of the symmetric group so that $\gamma hgh^{-1} = h \gamma g
h^{-1}$, we have the following corollary to Theorem 7.8.

\begin{cor}
Let $g' \in S_k$ be a $k$-cycle, let $g = (12 \cdots k)$ and let $h_1$
be the unique element of $S_k$ that fixes 1 and satisfies $g' =
h_1gh_1^{-1}$.  Then we have the following isomorphism of categories
\begin{eqnarray*}
T_{\gamma g'}^k = h_1 \circ T_{\gamma g}^k : \mathcal{ C}_w^{\gamma^k}(V) 
&\longrightarrow& \mathcal{ C}^{\gamma g'}_w(V^{\otimes k})\\
(M,Y) &\mapsto& (T_{\gamma g'}^k(M), Y_{\gamma g'}) = (M, Y_{h_1 \gamma 
gh_1^{-1}})
\end{eqnarray*}
where 
\[Y_{h_1 \gamma gh_1^{-1}} (v,z) = Y_{\gamma g} (h_1v,z) \]
and $Y_{\gamma g}$ is uniquely determined by (\ref{Y for gamma g}) and 
(\ref{Y for gamma g 2}).  Moreover, $T_{\gamma g'}^k$ preserves 
irreducible, admissible and ordinary objects.
\end{cor}

Finally we deal with the case $\gamma g$ when $g$ is an arbitrary
permutation in $S_k$.  As in Section 6 we write $g$ as a product of
disjoint cycles $g=hg'_1\cdots g'_ph^{-1}$ where the order of $g'_i$
is $k_i$ such that $\sum_{i}k_i=k$ and where $g_i$ permutes the
numbers $(\sum_{j = 1}^{i-1} k_j) + 1, (\sum_{j = 1}^{i-1} k_j) + 2, 
..., \sum_{j = 1}^i k_j$.  The following theorem generalizes Theorem 
\ref{t5.3}.

\begin{thm} Let $W_i$ be a weak $\gamma^{k_i}$-twisted 
$V$-module for $i=1,...,p$.  Given $g \in S_k$ and a decomposition
$g=hg'_1\cdots g'_ph^{-1}$ as above, $h \circ (T^{k_1}_{\gamma g_1} 
(W_1) \otimes \cdots \otimes T^{k_p}_{\gamma g_p}(W_p))$ is a weak 
$\gamma g$-twisted $\vtk$-module.  Moreover,

(1) If each $W_i$ is irreducible, then $h \circ (T^{k_1}_{\gamma g_1} 
(W_1) \otimes \cdots \otimes T^{k_p}_{\gamma g_p} (W_p))$ is irreducible.  

(2) If each $W_i$ is an admissible (resp., ordinary) $V$-module, 
then $h \circ (T^{k_1}_{\gamma g_1} (W_1) \otimes \cdots \otimes 
T^{k_p}_{\gamma g_p} (W_p))$ is an admissible (resp., ordinary) $\gamma
g$-twisted $\vtk$-module.

(3) If $V$ is a $\gamma^{k_i}$-rational vertex operator algebra for
$i=1,...,p$, then $\vtk$ is $\gamma g$-rational and all irreducible
$\gamma g$-twisted modules are given by $h \circ (T^{k_1}_{\gamma g_1} 
(W_1) \otimes \cdots \otimes T^{k_p}_{\gamma g_p} (W_p))$ where $W_i$ are
irreducible $\gamma^{k_i}$-twisted $V$-modules.
\end{thm}

Of course the analogous symmetries discussed in Remark 5.2 and Remark
6.5 hold for Corollary 7.9 and Theorem 7.10, respectively.

\end{document}